\documentclass[english]{article}
\usepackage[T1]{fontenc}
\usepackage[latin9]{inputenc}
\usepackage{geometry}
\usepackage{verbatim}
\usepackage{amsmath}
\usepackage{amssymb}
\usepackage{esint}

\makeatletter

\providecommand{\tabularnewline}{\\}

\newcommand{\lyxaddress}[1]{
\par {\raggedright #1
\vspace{1.4em}
\noindent\par}
}

\makeatother

\makeatother

\usepackage{babel}

\makeatother

\usepackage{babel}

\makeatother

\usepackage{babel}

\makeatother

\usepackage{babel}

\begin{document}

\title{$\,_{3}F_{4}$ hypergeometric functions as a sum of a product of
$\,_{2}F_{3}$ functions}

\author{Jack C. Straton}

\maketitle

\lyxaddress{Department of Physics, Portland State University, Portland, OR 97207-0751,
USA}

\lyxaddress{straton@pdx.edu}
\begin{abstract}
This paper shows that certain $\,_{3}F_{4}$ hypergeometric functions
may be expanded in sums of pair products of $\,_{2}F_{3}$ functions.
This expands the class of hypergeometric functions having summation
theorems beyond those expressible as pair-products of generalized
Whittaker functions, $\,_{2}F_{1}$ functions, and $\,_{3}F_{2}$ functions into the realm of  $\,_{P}F_{Q}$ functions where $P<Q$ for both  the summand and terms in the series.
In addition to its intrinsic value, this result has a specific application
in calculating the response of the atoms to laser stimulation in the
Strong Field Approximation.
\end{abstract}
\vspace{2pc}
 \\

\section{Introduction}

\noindent In 1942 Jackson\cite{Jackson} extended the concept of \emph{summation
theorems} over pair-products of functions such as Bessel functions
\cite{GR5 p. 992 No.  8.53} to sums of pair-products of the broader
class of Generalized Hypergeometric functions. Though his focus was
on hypergeometric functions of two variables (x,y), as a special case
he replaced $y=xq^{b}$, giving $\,_{2}F_{1}$ functions expanded
as pairs of $\,_{2}F_{1}$ functions (his eq. (I.55)) and $\,_{1}F_{1}$
functions (Whittaker functions) expanded as pairs of $\,_{1}F_{1}$
functions (his eq. (II.69)). In 1962 Ragab \cite{Ragab} found six
such expressions involving Slater's\cite{Slater} generalization of
Whittaker functions to $\,_{P}F_{P}$ functions having $P\geq1$,
all but one of which have $x^{2}$ rather than $x$ as the argument
in the summand, such as

\noindent \begin{eqnarray}
\,_{2}F_{3}\left(\frac{b}{2}+\frac{1}{2},\frac{b}{2};a+\frac{1}{2},\frac{c}{2}+\frac{1}{2},\frac{c}{2};\frac{x^{2}}{16}\right) &=&\sum_{r=0}^{\infty}\frac{(-1)^{r}x^{2r}(a)_{r}(b)_{r}(c-b)_{r}}{r!(2a)_{r}(c)_{2r}(c+r-1)_{r}} \nonumber \\ 
& \times & \,_{2}F_{2}(a+r,b+r;2a+r,c+2r;x)\,_{1}F_{1}\left(b+r;c+2r;-\frac{x}{2}\right)\end{eqnarray}


%
%
%

In addition to the summands that were $\,_{2}F_{3}$ hypergeometric
functions for many of his results, he expressed a $\,_{4}F_{5}$ function
as a sum of products of $\,_{2}F_{2}$ functions. Arun Verma\cite{Verma} rederived
Jackson's and some of Ragab's results in 1964 and added expansions
of $\,_{3}F_{2}$ functions as a product of a $\,_{2}F_{1}$ function
with another $\,_{3}F_{2}$ function. He also expressed a $\,_{5}F_{5}$
generalized Whittaker function as a sum of products of $\,_{2}F_{2}$
functions.

The present paper derives summation theorems for $\,_{3}F_{4}$ hypergeometric
functions  in terms of pair products
$\,_{2}F_{3}$  functions. While interesting in itself, this result has a specific
application in calculating the response of the atoms to laser stimulation
in the Strong Field Approximation (SFA) \cite{Reiss1,Reiss3,Reiss4,Reiss5,Reiss6,Faisel2}.
Whereas perturbation expansions will not converge if the applied laser
field is sufficiently large, the Strong Field Approximation (SFA)
is an analytical approximation that is non-perturbative. Keating, \cite{KeatingPhD}
applied it specifically to the production of the positive antihydrogen
ion.

\section{The transition amplitude}

Such SFA transition amplitudes involve integrals over differential
angles, in Keating's notation,

\begin{equation}
\Xi_{n}\left(k\alpha_{0},-\frac{z}{2}\right)=\int J_{n}^{2}\left(\mathbf{k}\cdot\boldsymbol{\alpha}_{0},-\frac{z}{2}\right)d\Omega\label{eq:xidef}\end{equation}
 of the Generalized Bessel function

\begin{equation}
J_{n}(x,y)=\sum_{h=-\infty}^{\infty}J_{n-2h}(x)J_{h}(y)\quad.\label{eq:jdef}\end{equation}

\subsection{Derivation of the Function $\Xi_{n}^{p}(k\alpha_{0},-\frac{z}{2})$}

We extend Keating's function somewhat with an additional cosine-squared
factor

\begin{equation}
\Xi_{n}^{p}\left(k\alpha_{0},-\frac{z}{2}\right)=\int\cos^{2p}\left(\theta\right)J_{n}^{2}\left(\mathbf{k}\cdot\boldsymbol{\alpha}_{0},-\frac{z}{2}\right)d\Omega\label{eq:xipdef}\end{equation}
 whose power \emph{2p} can be set to 0 to reproduce Keating's result
or retained with higher integers.

This integral can be solved analytically by expanding $\cos^{p}\left(\theta\right)J_{n}\left(\mathbf{k}\cdot\boldsymbol{\alpha}_{0},-\frac{z}{2}\right)$
in a Laplace series,\cite{Haber} that is, in terms of spherical
harmonics

\begin{equation}
\cos^{p}\left(\theta\right)J_{n}\left(\mathbf{k}\cdot\boldsymbol{\alpha}_{0},-\frac{z}{2}\right)=\sum_{lm}F_{lm}^{pn}(k\alpha_{0},-\frac{z}{2})Y_{lm}\left(\theta,\phi\right)\quad,\label{eq:Yexp}\end{equation}
 where, in a parallel notation to Keating's,

\begin{equation}
F_{lm}^{pn}(k\alpha_{0},-\frac{z}{2})=\int\cos^{p}\left(\theta\right)J_{n}\left(\mathbf{k}\cdot\boldsymbol{\alpha}_{0},-\frac{z}{2}\right)Y_{lm}^{*}\left(\theta,\phi\right)d\Omega.\label{eq:fint}\end{equation}

With these definitions,

\begin{equation}
\begin{array}{ccc}
\Xi_{n}^{p}\left(k\alpha_{0},-\frac{z}{2}\right) & = & \int\cos^{2p}\left(\theta\right)J_{n}^{2}\left(\mathbf{k}\cdot\boldsymbol{\alpha}_{0},-\frac{z}{2}\right)d\Omega\\
 & = & \int\sum_{lm}F_{lm}^{pn}(k\alpha_{0},-\frac{z}{2})Y_{lm}(\theta,\phi)\sum_{l'm'}F_{l'm'}^{pn}(k\alpha_{0},-\frac{z}{2})Y_{l'm'}^{*}\left(\theta,\phi\right)d\Omega\\
 & = & \sum_{lm}F_{lm}^{pn}(k\alpha_{0},-\frac{z}{2}))\sum_{l'm'}F_{l'm'}^{pn}(k\alpha_{0},-\frac{z}{2})\delta_{l'l}\delta_{m'm}\\
 & = & \sum_{lm}F_{lm}^{pn}(k\alpha_{0},-\frac{z}{2}))F_{lm}^{pn}(k\alpha_{0},-\frac{z}{2})\end{array}\label{eq:Chipn}\end{equation}
 Using the definition of the Generalized Bessel function, we find

\begin{equation}
F_{lm}^{pn}(k\alpha_{0},-\frac{z}{2})=\sum_{h=-\infty}^{\infty}\int\cos^{p}\left(\theta\right)J_{n-2h}(\mathbf{k}\cdot\boldsymbol{\alpha}_{0})Y_{lm}^{*}\left(\theta,\phi\right)d\Omega\, J_{h}(-\frac{z}{2})\quad.\label{eq:fhsum}\end{equation}

In the derivation leading to the $\,_{2}F_{3}$ functions that are
the focus of this work, we wish to avoid the infinities that come
from negative indices of the Bessel functions that we will integrate
over. So using $\delta=\left\{ \begin{array}{c}
1\\
0\end{array}\right\} $ for $n=\left\{ \begin{array}{c}
odd\\
even\end{array}\right\} $ and \cite{GR5 p. 979 No. 8.472}

\begin{equation}
J_{-N}\left(z\right)=\left(-1\right)^{N}J_{N}\left(z\right)\quad\left[\textrm{N is a natural number}\right]\label{eq:negJindex}\end{equation}
 this becomes

\begin{equation}
\begin{array}{ccc}
F_{lm}^{pn}\left(k\alpha_{0},-\frac{z}{2}\right) & = & \left({\displaystyle \sum_{h=-\infty}^{\left(n-\delta\right)/2}}+{\displaystyle \sum_{h=\left(n-\delta\right)/2+1}^{\infty}}\right)\int\cos^{p}\left(\theta\right)J_{n-2h}\left(\mathbf{k}\cdot\boldsymbol{\alpha}_{0}\right)Y_{lm}^{*}\left(\theta,\phi\right)d\Omega\, J_{h}\left(-\frac{z}{2}\right)\\
 & = & {\displaystyle \sum_{j=\left(n-\delta\right)/2}^{\infty}}\int\cos^{p}\left(\theta\right)J_{n+2j}\left(\mathbf{k}\cdot\boldsymbol{\alpha}_{0}\right)Y_{lm}^{*}\left(\theta,\phi\right)d\Omega\, J_{-j}\left(-\frac{z}{2}\right)\\
 & + & {\displaystyle \sum_{j=\left(n-\delta\right)/2+1}^{\infty}}\int\cos^{p}\left(\theta\right)J_{n-2j}\left(\mathbf{k}\cdot\boldsymbol{\alpha}_{0}\right)Y_{lm}^{*}\left(\theta,\phi\right)d\Omega\, J_{j}\left(-\frac{z}{2}\right)\\
 & = & {\displaystyle \sum_{j=\left(n-\delta\right)/2}^{\infty}}\int\cos^{p}\left(\theta\right)J_{2j+n}\left(\mathbf{k}\cdot\boldsymbol{\alpha}_{0}\right)Y_{lm}^{*}\left(\theta,\phi\right)d\Omega\, J_{-j}\left(-\frac{z}{2}\right)\\
 & + & \left(-1\right)^{n-2j}{\displaystyle \sum_{j=\left(n-\delta\right)/2+1}^{\infty}}\int\cos^{p}\left(\theta\right)J_{2j-n}\left(\mathbf{k}\cdot\boldsymbol{\alpha}_{0}\right)Y_{lm}^{*}\left(\theta,\phi\right)d\Omega\, J_{j}\left(-\frac{z}{2}\right)\end{array}\quad.\label{eq:fsplitsum}\end{equation}
 Then

\begin{equation}
\begin{array}{ccc}
F_{lm}^{pn}(k\alpha_{0},-\frac{z}{2}) & = & {\displaystyle \sum_{h=-\infty}^{\infty}}\mathbb{J}_{lm}^{p,n-2h}\left(k\alpha_{0}\right)J_{h}(-\frac{z}{2})\\
 & = & {\displaystyle \sum_{j=\left(n-\delta\right)/2}^{\infty}}\mathbb{J}_{lm}^{p,2j+n}\left(k\alpha_{0}\right)J_{-j}\left(-\frac{z}{2}\right)\\
 & + & \left(-1\right)^{n}{\displaystyle \sum_{j=\left(n-\delta\right)/2+1}^{\infty}}\mathbb{J}_{lm}^{p,2j-n}\left(k\alpha_{0}\right)\, J_{j}\left(-\frac{z}{2}\right)\end{array},\label{eq:fashollowj}\end{equation}
 where we have extended Keating's definition to

\begin{equation}
\mathbb{J}_{lm}^{p,2j\pm n}\left(k\alpha_{0}\right)=\int\cos^{p}\left(\theta\right)J_{2j\pm n}\left(\mathbf{k}\cdot\boldsymbol{\alpha}_{0}\right)Y_{lm}^{*}\left(\theta,\phi\right)d\Omega\quad.\label{eq:hollowjdef}\end{equation}
 We pause to note that if $p=0$, this is simply the Laplace series
of a conventional Bessel function with integer indices $J_{2j\pm n}\left(\mathbf{k}\cdot\boldsymbol{\alpha}_{0}\right)\equiv J_{2j\pm n}\left(k\alpha_{0}\cos\theta\right)$,
whose derivation we have not seen in the literature prior to Keating's.

Since this function is independent of the azimuthal angle, we can
reduce the Laplace series to a sum over Legendre polynomials. Using
the definition of the spherical harmonics

\begin{equation}
Y_{lm}^{*}\left(\theta,\phi\right)=(-1)^{l}\sqrt{\frac{2l+1}{4\pi}\frac{(l-m)!}{(l+m)!}}P_{l}(\cos\theta)e^{-im\phi}\quad,\label{eq:Ydef}\end{equation}
 we obtain

\begin{equation}
\mathbb{J}_{lm}^{p,2j\pm n}\left(k\alpha_{0}\right)=(-1)^{l}\sqrt{\frac{2l+1}{4\pi}\frac{(l-m)!}{(l+m)!}}\int\cos^{p}\left(\theta\right)J_{2j\pm n}\left(k\alpha_{0}\cos\theta\right)P_{l}\left(\cos\theta\right)d\left(\cos\theta\right)d\theta\int e^{-im\phi}d\phi\label{eq:hollowjints}\end{equation}
 or

\begin{equation}
\mathbb{J}_{lm}^{p,2j\pm n}\left(k\alpha_{0}\right)=(-1)^{l}\sqrt{\frac{2l+1}{4\pi}\frac{(l-m)!}{(l+m)!}}\mathcal{J}_{l}^{p,2j\pm n}(k\alpha_{0})\delta_{m0}\quad,\label{eq:curlytohollowj}\end{equation}
 where we again extend Keating's defined function to be

\begin{equation}
\mathcal{J}_{l}^{p,2j\pm n}\left(k\alpha_{0}\right)=2\pi\int_{-\pi}^{\pi}\cos^{p}\left(\theta\right)J_{2j\pm n}\left(k\alpha_{0}\cos\theta\right)P_{l}\left(\cos\theta\right)d\left(\cos\theta\right)\quad.\label{eq:curlyjints}\end{equation}
 Let $x=\cos\left(\theta\right)$, then

\begin{equation}
\mathcal{J}_{l}^{p,2j\pm n}(k\alpha_{0})=2\pi\int_{-1}^{1}x^{p}J_{2j\pm n}(k\alpha_{0}x)P_{l}(x)dx\quad.\label{eq:curleyjxint}\end{equation}

\section{The Fourier-Legendre series of a Bessel function of the first kind}

In a prior paper \cite{stra24a} we dropped the complicated indices
to cleanly express the Bessel function in a series of Legendre polynomials,
on the assumption \cite{Kellogg} that the series \begin{equation}
J_{N}\left(kx\right)=\sum_{L=0}^{\infty}a_{LN}\left(k\right)P_{L}(x)\quad\label{eq:Fourier-Legendre series}\end{equation}
converges uniformly,%
\footnote{Let D be a region in which the above series converges for each value
of x. Then the series can be said to converge uniformly in D if, for
every $\varepsilon>0$, there exists a number $N'(\varepsilon)$ such
that, for $n>N'$, it follows that

\[
\left.\left|J_{N}\left(kx\right)-\sum_{L=0}^{n}a_{LN}\left(k\right)P_{L}(x)\right|=\left|\sum_{L=n+1}^{\infty}a_{LN}\left(k\right)P_{L}(x)\right|\right.<\varepsilon\]

for all x in D. %
} where the coefficients are given by the orthogonality of the Legendre
polynomials,

\begin{equation}
a_{LN}\left(k\right)=\frac{2L+1}{2}\int_{-1}^{1}J_{N}\left(kx\right)P_{L}(x)dx\;.\label{eq:coefficients}\end{equation}

Following Keating's lead, but without the complicated indices, we
showed that 

\begin{equation}
\begin{array}{ccc}
a_{LN}\left(k\right) & = & \sqrt{\pi}(2L+1)2^{-L-1}i^{L-\text{N}}\sum_{M=0}^{\infty}\frac{\left(\left(-\frac{1}{4}\right)^{M}k^{L+2M}\right)}{2^{L+2M+1}\left(M!\Gamma\left(L+M+\frac{3}{2}\right)\right)}\\
 & \times & \left(1+(-1)^{L+2M+\text{N}}\right)\binom{L+2M}{\frac{1}{2}(L+2M-\text{N})}\\
 & = & \frac{\sqrt{\pi}2^{-2L-2}(2L+1)k^{L}i^{L-N}}{\Gamma\left(\frac{1}{2}(2L+3)\right)}\left(1+(-1)^{L+N}\right)\binom{L}{\frac{L-N}{2}}\\
 & \times & \,_{2}F_{3}\left(\frac{L}{2}+\frac{1}{2},\frac{L}{2}+1;L+\frac{3}{2},\frac{L}{2}-\frac{N}{2}+1,\frac{L}{2}+\frac{N}{2}+1;-\frac{k^{2}}{4}\right)\\
 & = & \sqrt{\pi}\,2^{-2L-2}(2L+1)k^{L}i^{L-N}\left(1+(-1)^{L+N}\right)\Gamma(L+1)\\
 & \times & \,_{2}\tilde{F}_{3}\left(\frac{L}{2}+\frac{1}{2},\frac{L}{2}+1;L+\frac{3}{2},\frac{L}{2}-\frac{N}{2}+1,\frac{L}{2}+\frac{N}{2}+1;-\frac{k^{2}}{4}\right)\end{array}\;,\label{eq:a_L_as_2F3}\end{equation}
where the final two steps are new with the prior work \cite{stra24a}.
Whenever \emph{$N>1$} is an integer larger than \emph{L}, and of
the same parity, the conventional form of the hypergeometric function
in the second form gives indeterminacies (ratios of infinities) in computation so we have
included the final form involving regularized hypergeometric functions
\cite{wolfram.com/07.26.26.0001.01}

\begin{equation}
\,_{2}F_{3}\left(a_{1},a_{2};b_{1},b_{2},b_{3};z\right)=\Gamma\left(b_{1}\right)\Gamma\left(b_{2}\right)\Gamma\left(b_{3}\right)\,_{2}\tilde{F}_{3}\left(a_{1},a_{2};b_{1},b_{2},b_{3};z\right)\label{eq:regularized-1}\end{equation}
and cancelled the $\Gamma\left(b_{i}\right)$ with gamma functions
in the denominators of the prefactors that each give infinities in
this case. 

For the special cases of $N=0$ and $1$ the order of the hypergeometric
functions is reduced since the parameters $a_{2}=b_{3}$ and $a_{1}=b_{2}$,
resp., giving

\begin{equation}
\begin{array}{ccc}
a_{L0}\left(k\right) & = & \frac{\sqrt{\pi}i^{L}2^{-2L-2}(2L+1)k^{L}}{\Gamma\left(\frac{1}{2}(2L+3)\right)}\left(1+(-1)^{L}\right)\binom{L}{\frac{L}{2}}\\
 & \times & \,_{1}F_{2}\left(\frac{L}{2}+\frac{1}{2};\frac{L}{2}+1,L+\frac{3}{2};-\frac{k^{2}}{4}\right)\\
 & = & \sqrt{\pi}i^{L}2^{-2L-2}(2L+1)k^{L}\Gamma\left(\frac{L}{2}+1\right)\left(1+(-1)^{L}\right)\binom{L}{\frac{L}{2}}\\
 & \times & \,_{1}\tilde{F}_{2}\left(\frac{L}{2}+\frac{1}{2};\frac{L}{2}+1,L+\frac{3}{2};-\frac{k^{2}}{4}\right)\end{array}\;,\label{eq:a_L0}\end{equation}

and

\begin{equation}
\begin{array}{ccc}
a_{L1}\left(k\right) & = & \frac{\sqrt{\pi}i^{L-1}2^{-2L-2}(2L+1)k^{L}}{\Gamma\left(\frac{1}{2}(2L+3)\right)}\left(1+(-1)^{L+1}\right)\binom{L}{\frac{L-1}{2}}\\
 & \times & \,_{1}F_{2}\left(\frac{L}{2}+1;\frac{L}{2}+\frac{3}{2},L+\frac{3}{2};-\frac{k^{2}}{4}\right)\\
 & = & i^{L-1}2^{-L-2}(2L+1)k^{L}\Gamma\left(\frac{L}{2}+1\right)\left(1+(-1)^{L+1}\right)\\
 & \times & \,_{1}\tilde{F}_{2}\left(\frac{L}{2}+1;\frac{L}{2}+\frac{3}{2},L+\frac{3}{2};-\frac{k^{2}}{4}\right)\end{array}\;,\label{eq:a_L1}\end{equation}
In each special case the first form involving a hypergeometric function
has no indeterminacies, but we include the regularized hypergeometric
function version for completeness.

Stepping outward in the string of definitions, with $k\rightarrow k\alpha_{0}$,

\begin{equation}
\begin{array}{ccc}
\mathbb{J}_{lm}^{0,2j\pm n}\left(k\alpha_{0}\right) & = & (-1)^{l}\sqrt{\frac{2l+1}{4\pi}\frac{(l-m)!}{(l+m)!}}\mathcal{J}_{l}^{0,2j\pm n}(k\alpha_{0})\delta_{m0},\\
 & = & \left(-1\right){}^{l}\sqrt{\frac{2l+1}{4\pi}}\frac{\pi^{3/2}2^{-2l}\left(k\alpha_{0}\right){}^{l}i^{(l\mp n)-2j}}{\Gamma\left(\frac{1}{2}(2l+3)\right)}\left((-1)^{2j+(l\pm n)}+1\right)\binom{l}{\frac{1}{2}((l\mp n)-2j)}\\
 & \times & \,_{2}F_{3}\left(\frac{l}{2}+\frac{1}{2},\frac{l}{2}+1;l+\frac{3}{2},-j+\left(\frac{l}{2}\mp\frac{n}{2}\right)+1,j+\left(\frac{l}{2}\pm\frac{n}{2}\right)+1;-\frac{1}{4}k^{2}\alpha_{0}^{2}\right)\end{array}\label{eq:hollowjas2F3}\end{equation}
 so that

\begin{equation}
\begin{array}{ccc}
F_{lm}^{0n}(k\alpha_{0},-\frac{z}{2}) & = & {\displaystyle \sum_{j=\left(n-\delta\right)/2}^{\infty}\mathbb{J}_{lm}^{0,2j+n}\left(k\alpha_{0}\right)}J_{-j}\left(-\frac{z}{2}\right)\\
 & + & \left(-1\right)^{n}{\displaystyle \sum_{j=\left(n-\delta\right)/2+1}^{\infty}}\mathbb{J}_{lm}^{0,2j-n}\left(k\alpha_{0}\right)\, J_{j}\left(-\frac{z}{2}\right)\\
 & = & {\displaystyle \sum_{j=\left(n-\delta\right)/2}^{\infty}}\left(-1\right){}^{l}\sqrt{\frac{2l+1}{4\pi}}\frac{\pi^{3/2}2^{-2l}\left(k\alpha_{0}\right){}^{l}i^{(l-n)-2j}}{\Gamma\left(\frac{1}{2}(2l+3)\right)}\left((-1)^{2j+(l+n)}+1\right)\binom{l}{\frac{1}{2}((l-n)-2j)}J_{-j}\left(-\frac{z}{2}\right)\\
 & \times & \,_{2}F_{3}\left(\frac{l}{2}+\frac{1}{2},\frac{l}{2}+1;l+\frac{3}{2},-j+\left(\frac{l}{2}-\frac{n}{2}\right)+1,j+\left(\frac{l}{2}+\frac{n}{2}\right)+1;-\frac{1}{4}k^{2}\alpha_{0}^{2}\right)\\
 & + & \left(-1\right)^{n}{\displaystyle \sum_{j=\left(n-\delta\right)/2+1}^{\infty}}\left(-1\right){}^{l}\sqrt{\frac{2l+1}{4\pi}}\frac{\pi^{3/2}2^{-2l}\left(k\alpha_{0}\right){}^{l}i^{(l+n)-2j}}{\Gamma\left(\frac{1}{2}(2l+3)\right)}\left((-1)^{2j+(l-n)}+1\right)\binom{l}{\frac{1}{2}((l+n)-2j)} \\
 & \times & \, J_{j}\left(-\frac{z}{2}\right) \,_{2}F_{3}\left(\frac{l}{2}+\frac{1}{2},\frac{l}{2}+1;l+\frac{3}{2},-j+\left(\frac{l}{2}+\frac{n}{2}\right)+1,j+\left(\frac{l}{2}-\frac{n}{2}\right)+1;-\frac{1}{4}k^{2}\alpha_{0}^{2}\right)\end{array},\label{eq:Fas2F3}\end{equation}
 and finally

\begin{equation}
\begin{array}{ccc}
\Xi_{n}^{0}\left(k\alpha_{0},-\frac{z}{2}\right) & = & \sum_{l0}F_{l0}^{0n}(k\alpha_{0},-\frac{z}{2})F_{l0}^{0n}(k\alpha_{0},-\frac{z}{2})\\
 & = & \sum_{l0}\frac{2l+1}{4\pi}\left({\displaystyle \sum_{j=-\left(n-\delta\right)/2}^{\infty}\quad}{\displaystyle \sum_{M=-\left(n-\delta\right)/2}^{\infty}}J_{-j}\left(-\frac{z}{2}\right)J_{-M}\left(-\frac{z}{2}\right)\right.\\
 & \times & \frac{\pi^{3/2}2^{-2l}\left(k\alpha_{0}\right){}^{l}i^{(l-n)-2j}}{\Gamma\left(\frac{1}{2}(2l+3)\right)}\left((-1)^{2j+(l+n)}+1\right)\binom{l}{\frac{1}{2}((l-n)-2j)}\\
 & \times & \frac{\pi^{3/2}2^{-2l}\left(k\alpha_{0}\right){}^{l}i^{l-2M-n}}{\Gamma\left(\frac{1}{2}(2l+3)\right)}\left((-1)^{l+2M+n}+1\right)\binom{l}{\frac{1}{2}(l-2M-n)}\\
 & \times & \,_{2}F_{3}\left(\frac{l}{2}+\frac{1}{2},\frac{l}{2}+1;l+\frac{3}{2},-j+\left(\frac{l}{2}-\frac{n}{2}\right)+1,j+\left(\frac{l}{2}+\frac{n}{2}\right)+1;-\frac{1}{4}k^{2}\alpha_{0}^{2}\right)\\
 & \times & _{2}F_{3}\left(\frac{l}{2}+\frac{1}{2},\frac{l}{2}+1;l+\frac{3}{2},\frac{l}{2}-M-\frac{n}{2}+1,\frac{l}{2}+M+\frac{n}{2}+1;-\frac{1}{4}k^{2}\alpha_{0}^{2}\right)\\
 & + & \left(-1\right)^{n}{\displaystyle \sum_{j=\left(n-\delta\right)/2+1}^{\infty}}{\displaystyle \sum_{M=-\left(n-\delta\right)/2}^{\infty}}J_{j}\left(-\frac{z}{2}\right)J_{-M}\left(-\frac{z}{2}\right)\\
 & \times & \frac{\pi^{3/2}2^{-2l}\left(k\alpha_{0}\right){}^{l}i^{(l-n)-2j}}{\Gamma\left(\frac{1}{2}(2l+3)\right)}\left((-1)^{2j+(l+n)}+1\right)\binom{l}{\frac{1}{2}((l-n)-2j)}\\
 & \times & \frac{\pi^{3/2}2^{-2l}\left(k\alpha_{0}\right){}^{l}i^{l-2M+n}}{\Gamma\left(\frac{1}{2}(2l+3)\right)}\left((-1)^{l+2M-n}+1\right)\binom{l}{\frac{1}{2}(l-2M+n)}\\
 & \times & \,_{2}F_{3}\left(\frac{l}{2}+\frac{1}{2},\frac{l}{2}+1;l+\frac{3}{2},-j+\left(\frac{l}{2}-\frac{n}{2}\right)+1,j+\left(\frac{l}{2}+\frac{n}{2}\right)+1;-\frac{1}{4}k^{2}\alpha_{0}^{2}\right)\\
 & \times & \,_{2}F_{3}\left(\frac{l}{2}+\frac{1}{2},\frac{l}{2}+1;l+\frac{3}{2},\frac{l}{2}+M-\frac{n}{2}+1,\frac{l}{2}-M+\frac{n}{2}+1;-\frac{1}{4}k^{2}\alpha_{0}^{2}\right)\quad.\\
 & + & \left(-1\right)^{n}{\displaystyle \sum_{j=-\left(n-\delta\right)/2}^{\infty}}{\displaystyle \sum_{M=\left(n-\delta\right)/2+1}^{\infty}}\left(-\frac{z}{2}\right)J_{M}\left(-\frac{z}{2}\right)\\
 & \times & \frac{\pi^{3/2}2^{-2l}\left(k\alpha_{0}\right){}^{l}i^{(l+n)-2j}}{\Gamma\left(\frac{1}{2}(2l+3)\right)}\left((-1)^{2j+(l-n)}+1\right)\binom{l}{\frac{1}{2}((l+n)-2j)}\\
 & \times & \frac{\pi^{3/2}2^{-2l}\left(k\alpha_{0}\right){}^{l}i^{l-2M-n}}{\Gamma\left(\frac{1}{2}(2l+3)\right)}\left((-1)^{l+2M+n}+1\right)\binom{l}{\frac{1}{2}(l-2M-n)}\\
 & \times & \,_{2}F_{3}\left(\frac{l}{2}+\frac{1}{2},\frac{l}{2}+1;l+\frac{3}{2},-j+\left(\frac{l}{2}+\frac{n}{2}\right)+1,j+\left(\frac{l}{2}-\frac{n}{2}\right)+1;-\frac{1}{4}k^{2}\alpha_{0}^{2}\right)\\
 & \times & _{2}F_{3}\left(\frac{l}{2}+\frac{1}{2},\frac{l}{2}+1;l+\frac{3}{2},\frac{l}{2}-M-\frac{n}{2}+1,\frac{l}{2}+M+\frac{n}{2}+1;-\frac{1}{4}k^{2}\alpha_{0}^{2}\right)\\
 & + & {\displaystyle \sum_{j=\left(n-\delta\right)/2+1}^{\infty}}{\displaystyle {\displaystyle \sum_{M=\left(n-\delta\right)/2+1}^{\infty}}}J_{j}\left(-\frac{z}{2}\right)J_{M}\left(-\frac{z}{2}\right)\\
 & \times & \frac{\pi^{3/2}2^{-2l}\left(k\alpha_{0}\right){}^{l}i^{(l-n)-2j}}{\Gamma\left(\frac{1}{2}(2l+3)\right)}\left((-1)^{2j+(l+n)}+1\right)\binom{l}{\frac{1}{2}((l-n)-2j)}\\
 & \times & \frac{\pi^{3/2}2^{-2l}\left(k\alpha_{0}\right){}^{l}i^{l-2M-n}}{\Gamma\left(\frac{1}{2}(2l+3)\right)}\left((-1)^{l+2M+n}+1\right)\binom{l}{\frac{1}{2}(l-2M-n)}\\
 & \times & \,_{2}F_{3}\left(\frac{l}{2}+\frac{1}{2},\frac{l}{2}+1;l+\frac{3}{2},-j+\left(\frac{l}{2}-\frac{n}{2}\right)+1,j+\left(\frac{l}{2}+\frac{n}{2}\right)+1;-\frac{1}{4}k^{2}\alpha_{0}^{2}\right)\\
 & \times & \left._{2}F_{3}\left(\frac{l}{2}+\frac{1}{2},\frac{l}{2}+1;l+\frac{3}{2},\frac{l}{2}-M-\frac{n}{2}+1,\frac{l}{2}+M+\frac{n}{2}+1;-\frac{1}{4}k^{2}\alpha_{0}^{2}\right)\right)\end{array}\label{eq:chi0sum}\end{equation}

\section{A more direct approach}

On the other hand, we can transform from Bessel to hypergeometric
functions using \cite{PBM3 p. 594 No. 7.13.1.1,Luke1 p. 212 No. 6.2.7.1,wolfram.com/03.01.26.0002.01}

\begin{equation}
J_{\nu}(z)=\frac{\left(\frac{z}{2}\right)^{\nu}\,_{0}F_{1}\left(;\nu+1;-\frac{z^{2}}{4}\right)}{\Gamma(\nu+1)}\label{eq:Jas0F1}\end{equation}
 and combine pairs\cite{EH p. 185 No. 4.3.2,wolfram.com/07.17.26.0004.01,Luke1  p. 228 No. 6.4.1.261}
 
 \begin{equation}
\,_{0}F_{1}(;b;z)\,_{0}F_{1}(;c;z)=\,_{2}F_{3}\left(\frac{b}{2}+\frac{c}{2}-\frac{1}{2},\frac{b}{2}+\frac{c}{2};b,c,b+c-1;4z\right)\label{eq:0F1 0F1 as 2F3}\end{equation}
 so that\cite{Luke1  p. 216 No. 6.2.7.39} \begin{equation}
J_{\mu}(z)J_{\nu}(z)=\frac{2^{-\mu-\nu}z^{\mu+\nu}}{\Gamma(\mu+1)\Gamma(\nu+1)}\,_{2}F_{3}\left(\frac{\mu}{2}+\frac{\nu}{2}+\frac{1}{2},\frac{\mu}{2}+\frac{\nu}{2}+1;\mu+1,\nu+1,\mu+\nu+1;-z^{2}\right)\quad.\label{eq:JJ as 2F3}\end{equation}

Then

\begin{equation}
\begin{array}{ccc}
J_{n}^{2} & & \hspace{-1.2cm} \left(k\text{\ensuremath{\alpha_{0}}}\cos(\theta),-\frac{z}{2}\right)  =  \left({\displaystyle \sum_{j=-\left(n-\delta\right)/2}^{\infty}}{\displaystyle \sum_{M=-\left(n-\delta\right)/2}^{\infty}}J_{n+2j}\left(k\text{\ensuremath{\alpha_{0}}}\cos(\theta)\right)J_{n+2M}(k\text{\ensuremath{\alpha_{0}}}\cos(\theta))  J_{-j}\left(-\frac{z}{2}\right)J_{-M}\left(-\frac{z}{2}\right) \right. \\
 & + & \left.\left(-1\right)^{n}{\displaystyle \sum_{j=\left(n-\delta\right)/2+1}^{\infty}}{\displaystyle \sum_{M=-\left(n-\delta\right)/2}^{\infty}}J_{2j-n}\left(k\text{\ensuremath{\alpha_{0}}}\cos(\theta)\right)J_{n+2M}(k\text{\ensuremath{\alpha_{0}}}\cos(\theta))J_{j}\left(-\frac{z}{2}\right)J_{-M}\left(-\frac{z}{2}\right)\right)\\
 & + & \left(\left(-1\right)^{n}{\displaystyle \sum_{j=-\left(n-\delta\right)/2}^{\infty}}{\displaystyle \sum_{M=\left(n-\delta\right)/2+1}^{\infty}}J_{n+2j}\left(k\text{\ensuremath{\alpha_{0}}}\cos(\theta)\right)J_{2M-n}(k\text{\ensuremath{\alpha_{0}}}\cos(\theta))J_{-j}\left(-\frac{z}{2}\right)J_{M}\left(-\frac{z}{2}\right)\right.\\
 & + & \left.{\displaystyle \sum_{j=\left(n-\delta\right)/2+1}^{\infty}}{\displaystyle {\displaystyle \sum_{M=\left(n-\delta\right)/2+1}^{\infty}}}J_{2j-n}\left(k\text{\ensuremath{\alpha_{0}}}\cos(\theta)\right)J_{2M-n}(k\text{\ensuremath{\alpha_{0}}}\cos(\theta))J_{j}\left(-\frac{z}{2}\right)J_{M}\left(-\frac{z}{2}\right)\right)\\
\\\\ & = & {\displaystyle \sum_{j=-\left(n-\delta\right)/2}^{\infty}}{\displaystyle \sum_{M=-\left(n-\delta\right)/2}^{\infty}}\frac{2^{-2j-2M-2n}(k\alpha_{0}\cos(\theta))^{2j+2M+2n}}{\Gamma(2j+n+1)\Gamma(2M+n+1)}J_{-j}\left(-\frac{z}{2}\right)J_{-M}\left(-\frac{z}{2}\right)\\
 & \times & \,_{2}F_{3}\left(j+M+n+\frac{1}{2},j+M+n+1;2j+n+1,2M+n+1,2j+2M+2n+1;-k^{2}\text{\ensuremath{\alpha_{0}}}^{2}\cos^{2}(\theta)\right)\\
 & + & \left(-1\right)^{n}{\displaystyle \sum_{j=\left(n-\delta\right)/2+1}^{\infty}}{\displaystyle \sum_{M=-\left(n-\delta\right)/2}^{\infty}}\frac{2^{-2j-2M}\left(k\alpha_{0}\cos(\theta)\right){}^{2j+2M}}{\Gamma(2j-n+1)\Gamma(2M+n+1)}J_{j}\left(-\frac{z}{2}\right)J_{-M}\left(-\frac{z}{2}\right)\\
 & \times & \,_{2}F_{3}\left(j+M+\frac{1}{2},j+M+1;2j+2M+1,2j-n+1,2M+n+1;-k^{2}\cos^{2}(\theta)\alpha_{0}^{2}\right)\\
 & + & \left(-1\right)^{n}{\displaystyle \sum_{j=-\left(n-\delta\right)/2}^{\infty}}{\displaystyle \sum_{M=\left(n-\delta\right)/2+1}^{\infty}}\frac{2^{-2j-2M}\left(k\alpha_{0}\cos(\theta)\right){}^{2j+2M}}{\Gamma(2j+n+1)\Gamma(2M-n+1)}J_{-j}\left(-\frac{z}{2}\right)J_{M}\left(-\frac{z}{2}\right)\\
 & \times & \,_{2}F_{3}\left(j+M+\frac{1}{2},j+M+1;2j+2M+1,2M-n+1,2j+n+1;-k^{2}\cos^{2}(\theta)\alpha_{0}^{2}\right)\\
 & + & {\displaystyle \sum_{j=\left(n-\delta\right)/2+1}^{\infty}}{\displaystyle {\displaystyle \sum_{M=\left(n-\delta\right)/2+1}^{\infty}}}\frac{2^{-2j-2M+2n}\left(k\alpha_{0}\cos(\theta)\right){}^{2j+2M-2n}}{\Gamma(2j-n+1)\Gamma(2M-n+1)}J_{j}\left(-\frac{z}{2}\right)J_{M}\left(-\frac{z}{2}\right)\\
 & \times & \,_{2}F_{3}\left(j+M-n+\frac{1}{2},j+M-n+1;2j+2M-2n+1,2j-n+1,2M-n+1;-k^{2}\cos^{2}(\theta)\alpha_{0}^{2}\right)\end{array}\quad.\label{eq:genJdefF3}\end{equation}

The integral we wish to do is

\begin{eqnarray}
\Xi_{n}^{p}\left(k\alpha_{0},-\frac{z}{2}\right) & = & \int\cos^{2p}\left(\theta\right)J_{n}^{2}\left(\mathbf{k}\cdot\boldsymbol{\alpha}_{0},-\frac{z}{2}\right)d\Omega=2\pi\int_{-\pi}^{\pi}\cos^{2p}\left(\theta\right)J_{n}^{2}\left(k\alpha_{0}\cos(\theta),-\frac{z}{2}\right)d\left(\cos(\theta)\right)\label{eq:angintsetup}\\
 &  & =2\pi\int_{-1}^{1}f\left(u^{2}\right)du=4\pi\int_{0}^{1}f\left(u^{2}\right)du=2\pi\int_{0}^{1}f\left(y\right)y^{-1/2}dy\quad.\nonumber 
 \end{eqnarray}
 We can then use \cite{PBM3 p. 334 No. 2.22.2.1}
 
 \begin{align}
\int_{0}^{a}y^{\alpha-1} & (a-y)^{\beta-1}\,_{p}F_{q}(a_{1},\ldots,\, a_{p},b_{1},\ldots,\, b_{q};-\omega y)\, dy & =\label{eq:PBM3 p. 334 No. 2.22.2.1}\\
& =  \frac{\Gamma(\alpha)\Gamma(\beta)a^{\alpha+\beta-1}}{\Gamma(\alpha+\beta)}\,_{p+1}F_{q+1}(a_{1},\ldots,\, a_{p},\alpha;b_{1},\ldots,\, b_{q},\alpha+\beta;-a\omega)\quad\left[\Re(\alpha)>0\land\Re(\beta)>0\land a>0\right]\end{align}
 %
with  $\alpha=\left\{ j+M+n+1,\: j+M+1,\: j+M+1,\: j+M-n+1\right\} +p$ for  the four terms and with $a=1$ and $\beta=1$ for each of the four.

Then

\begin{eqnarray}
\Xi_{n} & & \hspace{-0.9cm}  \left(k\alpha_{0},-\frac{z}{2}\right)  =  \int J_{n}^{2}\left(\mathbf{k}\cdot\boldsymbol{\alpha}_{0},-\frac{z}{2}\right)d\Omega\nonumber \\
 & = & 2\pi{\displaystyle \sum_{j=-\left(n-\delta\right)/2}^{\infty}}{\displaystyle \sum_{M=-\left(n-\delta\right)/2}^{\infty}}\frac{2^{-2j-2M-2n+1}\left(k\alpha_{0}\right){}^{2j+2M+2n}}{(2j+2M+2n+1)\Gamma(2j+n+1)\Gamma(2M+n+1)}J_{-j}\left(-\frac{z}{2}\right)J_{-M}\left(-\frac{z}{2}\right) \nonumber \\
 & \times & \,_{3}F_{4}\left(j+M+n+\frac{1}{2},j+M+n+\frac{1}{2},j+M+n+1;    \right.  \nonumber \\ 
 & & \left.  \hspace{3.1cm} 2j+n+1,j+M+n+\frac{3}{2},2M+n+1,2j+2M+2n+1;-k^{2}\alpha_{0}^{2}\right)\nonumber \\
 & + & 2\pi\left(-1\right)^{n}{\displaystyle \sum_{j=\left(n-\delta\right)/2+1}^{\infty}}{\displaystyle \sum_{M=-\left(n-\delta\right)/2}^{\infty}}\frac{2^{-2j-2M+1}\left(k\alpha_{0}\right){}^{2j+2M}}{(2j+2M+1)\Gamma(2j-n+1)\Gamma(2M+n+1)}J_{j}\left(-\frac{z}{2}\right)J_{-M}\left(-\frac{z}{2}\right)\nonumber \\
 & \times & \,_{3}F_{4}\left(j+M+\frac{1}{2},j+M+\frac{1}{2},j+M+1;j+M+\frac{3}{2},2j+2M+1,2j-n+1,2M+n+1;-k^{2}\alpha_{0}^{2}\right)\nonumber \\
 & + & 2\pi\left(-1\right)^{n}{\displaystyle \sum_{j=-\left(n-\delta\right)/2}^{\infty}}{\displaystyle \sum_{M=\left(n-\delta\right)/2+1}^{\infty}}\frac{2^{-2j-2M+1}\left(k\alpha_{0}\right){}^{2j+2M}}{(2j+2M+1)\Gamma(2j+n+1)\Gamma(2M-n+1)}J_{-j}\left(-\frac{z}{2}\right)J_{M}\left(-\frac{z}{2}\right)\nonumber \\
 & \times & \,_{3}F_{4}\left(j+M+\frac{1}{2},j+M+\frac{1}{2},j+M+1;j+M+\frac{3}{2},2j+2M+1,2M-n+1,2j+n+1;-k^{2}\alpha_{0}^{2}\right)\nonumber \\
 & + & 2\pi{\displaystyle \sum_{j=\left(n-\delta\right)/2+1}^{\infty}}{\displaystyle {\displaystyle \sum_{M=\left(n-\delta\right)/2+1}^{\infty}}}\frac{2^{-2j-2M+2n+1}\left(k\alpha_{0}\right){}^{2j+2M-2n}}{(2j+2M-2n+1)\Gamma(2j-n+1)\Gamma(2M-n+1)}J_{j}\left(-\frac{z}{2}\right)J_{M}\left(-\frac{z}{2}\right)\nonumber \\
 & \times & \,_{3}F_{4}\left(j+M-n+\frac{1}{2},j+M-n+\frac{1}{2},j+M-n+1; \right.  \nonumber \\ 
 & & \left.  \hspace{3.1cm}  2j+2M-2n+1,2j-n+1,j+M-n+\frac{3}{2},2M-n+1;-k^{2}\alpha_{0}^{2}\right)\label{eq:chifinal}\end{eqnarray}

Examination of equations (\ref{eq:chi0sum}) and (\ref{eq:chifinal})
shows that

\begin{equation}
\begin{array}{ccc}
\,_{3}F_{4} & & \hspace{-1.0cm}  \left(\frac{\mu}{2}+\frac{\nu}{2}+\frac{1}{2},\frac{\mu}{2}+\frac{\nu}{2}+\frac{1}{2},\frac{\mu}{2}+\frac{\nu}{2}+1;\mu+1,\frac{\mu}{2}+\frac{\nu}{2}+\frac{3}{2},\nu+1,\mu+\nu+1;z\right)  =  (\mu+\nu+1)\Gamma(\mu+1)\Gamma(\nu+1)\\
 & \times & {\displaystyle \sum_{L=0}^{\infty}\frac{\pi^{2}(2L+1)\left((-1)^{L+\mu}+1\right)\left((-1)^{L+\nu}+1\right)\Gamma(2L+2)^{2}i^{4L-2\mu-2\nu}2^{-8L+\mu+\nu-6}}{\Gamma\left(L+\frac{3}{2}\right)^{2}\Gamma\left(\frac{1}{2}(2L+3)\right)^{2}\Gamma\left(\frac{L}{2}-\frac{\mu}{2}+1\right)\Gamma\left(\frac{L}{2}+\frac{\mu}{2}+1\right)\Gamma\left(\frac{L}{2}-\frac{\nu}{2}+1\right)\Gamma\left(\frac{L}{2}+\frac{\nu}{2}+1\right)} } \\
 & \times &  z^{L-(\mu+\nu)/2} \,_{2}F_{3}\left(\frac{L}{2}+\frac{1}{2},\frac{L}{2}+1;L+\frac{3}{2},\frac{L}{2}-\frac{\mu}{2}+1,\frac{L}{2}+\frac{\mu}{2}+1;\frac{z}{4}\right) \\
 & \times &\,_{2}F_{3}\left(\frac{L}{2}+\frac{1}{2},\frac{L}{2}+1;L+\frac{3}{2},\frac{L}{2}-\frac{\nu}{2}+1,\frac{L}{2}+\frac{\nu}{2}+1;\frac{z}{4}\right)\\
 & = & (\mu+\nu+1)\Gamma(\mu+1)\Gamma(\nu+1)\\
 & \times & {\displaystyle \sum_{L=0}^{\infty}\frac{\pi^{2}(2L+1)\left((-1)^{L+\mu}+1\right)\left((-1)^{L+\nu}+1\right)\Gamma(2L+2)^{2}i^{4L-2\mu-2\nu}2^{-8L+\mu+\nu-6}}{\Gamma\left(\frac{1}{2}(2L+3)\right)^{2}}}  \\
 & \times & z^{L-(\mu+\nu)/2} \,_{2}\tilde{F}_{3}\left(\frac{L}{2}+\frac{1}{2},\frac{L}{2}+1;L+\frac{3}{2},\frac{L}{2}-\frac{\mu}{2}+1,\frac{L}{2}+\frac{\mu}{2}+1;\frac{z}{4}\right) \\
 & \times &\,_{2}\tilde{F}_{3}\left(\frac{L}{2}+\frac{1}{2},\frac{L}{2}+1;L+\frac{3}{2},\frac{L}{2}-\frac{\nu}{2}+1,\frac{L}{2}+\frac{\nu}{2}+1;\frac{z}{4}\right)\end{array}\quad.\label{eq:final0}\end{equation}
Numerical checks show that the right-hand side requires as few as
three nonzero terms -- $L=0,\,2,\,4$ or $L=1,\,3,\,5$, depending
on the parity of $\mu$ -- in the sum to get seven-digit accuracy
for $\mu+\nu<10$,

\begin{tabular}{|l|l|c|c|c|}
\hline 
Left-hand side of   (\ref{eq:final0})  & Right-hand side  of   (\ref{eq:final0}) & $\mu$ & $\nu$ & z\tabularnewline
\hline
\hline 
1.028881345119003 & 1.028881345119001 & 0 & 0 & 0.17\tabularnewline
\hline 
1.0344878191148 & 1.0344878191146 & 0 & 2 & 0.17\tabularnewline
\hline 
1.0369001971 & 1.0369001970 & 0 & 4 & 0.17\tabularnewline
\hline 
1.020434382759 & 1.020434382749 & 2 & 2 & 0.17\tabularnewline
\hline 
1.01777403 & 1.01777403 & 2 & 4 & 0.17\tabularnewline
\hline 
1.0140011 & 1.0140009 & 4 & 4 & 0.17\tabularnewline
\hline 
1.0258250454427744 & 1.0258250454427744 & 1 & 1 & 0.17\tabularnewline
\hline 
1.0230034607369 & 1.0230034607370 & 1 & 3 & 0.17\tabularnewline
\hline 
1.022243424630 & 1.022243424628 & 1 & 5 & 0.17\tabularnewline
\hline 
1.016657722535 & 1.016657722534 & 3 & 3 & 0.17\tabularnewline
\hline 
1.014587307 & 1.014587305 & 3 & 5 & 0.17\tabularnewline
\hline 
1.01205576 & 1.01205571 & 5 & 5 & 0.17\tabularnewline
\hline 
23.049 & 23.044 & 0 & 0 & 17.0\tabularnewline
\hline 
1.00013910008 & 1.00013910005 & 4 & 4 & 0.0017\tabularnewline
\hline
\end{tabular}

Table 1. The left and right sides of (\ref{eq:final0}) when the variable
is set to the arbitrary value $z=0.17$ and we include only three
nonzero terms in the sum, shown through the digit with which the two
sides disagree. When $z=.0017$ the accuracy increases, as seen in
the last line, while for $z=17$ the accuracy decreases; the penultimate
line. The results are symmetrical with respect to $\mu$ and $\nu$.

\vspace{0.2cm}

Parameters $\mu$ and $\nu$ must be of the same parity or the sum
is zero. They also must be non-negative to avoid infinities. The conventional
hypergeometric functions on the right-hand side give indeterminacies
in computation (infinities divided by infinities) unless $\mu$ and $\nu$ are both zero or one, so we
relied on the second version, using regularized hypergeometric functions
-- having cancelled the $\Gamma\left(b_{i}\right)$ with gamma functions
in the denominators of the prefactors that also give infinities in
this case -- for numerical checks. 

The derivation of the Fourier-Legendre series for $J_{\mu}(z)$ in
the prior paper \cite{stra24a}, upon which the present work relies,
was restricted to integer $\mu$ (and, hence, $\nu$ in the present
work) at two places, equations (3) and (10) of that paper, with the
latter reproduced below in (\ref{eq:GR5 p. 417 No. 3.631.17}). Equation
(3) of that paper is an integer-restricted version (the first term)
of the more general integral representation due to Heine of the Bessel
function \cite{MagnusOberhettingerSoni} 

\begin{equation}
\pi J_{\mu}(kx\cos(\theta))=e^{i\mu\pi/2}\left[\int_{0}^{\pi}dt\, e^{-ikx\cos(\theta)\cos(t)}\cos(\mu t)-\sin(\mu\pi)\int_{0}^{\infty}dt\, e^{ikx\cos(\theta)\cosh(t)-\mu t}\right]\label{eq:heine_general}\end{equation}
and, hence, is not the blockage to generalization. The essential blockage
is that we found no non-integer version of (\ref{eq:GR5 p. 417 No. 3.631.17}),
but perhaps one can be found. If found, one would have to investigate
whether the second term of (\ref{eq:heine_general}) can be integrated
over both $\cos(\theta)$ and $t$ if one wanted to remove the  integer restriction on $\mu$ and $\nu$.

When $\mu=\nu=0$ the order of the hypergeometric functions is reduced
since the parameters $a_{3}=b_{4}$ on the left-hand side and $a_{2}=b_{3}$
on the right-hand side, giving the special case,

\begin{equation}
\begin{array}{ccc}
\,_{2}F_{3}\left(\frac{1}{2},\frac{1}{2};1,1,\frac{3}{2};z\right) & = & \sum_{L=0}^{\infty}\frac{\pi^{2}i^{4L}2^{-8L-6}\left((-1)^{L}+1\right)^{2}(2L+1)\Gamma(2L+2)^{2}}{\Gamma\left(\frac{L}{2}+1\right)^{4}\Gamma\left(L+\frac{3}{2}\right)^{2}\Gamma\left(\frac{1}{2}(2L+3)\right)^{2}}z^{L}\,_{1}F_{2}\left(\frac{L}{2}+\frac{1}{2};\frac{L}{2}+1,L+\frac{3}{2};\frac{z}{4}\right)^{2}\\
 & = & \sum_{L=0}^{\infty}\frac{\pi^{2}i^{4L}2^{-8L-6}\left((-1)^{L}+1\right)^{2}(2L+1)\Gamma(2L+2)^{2}}{\Gamma\left(\frac{L}{2}+1\right)^{2}\Gamma\left(\frac{1}{2}(2L+3)\right)^{2}}z^{L}\,_{1}\tilde{F}_{2}\left(\frac{L}{2}+\frac{1}{2};L+\frac{3}{2},\frac{L}{2}+1;\frac{z}{4}\right)^{2}\end{array}\quad.\label{eq:final0_00}\end{equation}
When $\mu=\nu=1$ the order of the hypergeometric functions is also
reduced since the parameters $a_{3}=b_{3}$ on the left-hand side
and $a_{2}=b_{2}$ on the right-hand side, giving the special case,

\begin{equation}
\begin{array}{ccc}
\,_{2}F_{3}\left(\frac{3}{2},\frac{3}{2};2,\frac{5}{2},3;z\right) & = & 3\:\sum_{L=0}^{\infty}\frac{\pi^{2}i^{4L-4}2^{-8L-4}\left((-1)^{L+1}+1\right)^{2}(2L+1)\Gamma(2L+2)^{2}}{\Gamma\left(\frac{L}{2}+\frac{1}{2}\right)^{2}\Gamma\left(\frac{L}{2}+\frac{3}{2}\right)^{2}\Gamma\left(L+\frac{3}{2}\right)^{2}\Gamma\left(\frac{1}{2}(2L+3)\right)^{2}}z^{L-1}\,_{1}F_{2}\left(\frac{L}{2}+1;\frac{L}{2}+\frac{3}{2},L+\frac{3}{2};\frac{z}{4}\right)^{2}\\
 & = & 3\:\sum_{L=0}^{\infty}\frac{\pi^{2}i^{4L-4}2^{-8L-4}\left((-1)^{L+1}+1\right)^{2}(2L+1)\Gamma(2L+2)^{2}}{\Gamma\left(\frac{L}{2}+\frac{1}{2}\right)^{2}\Gamma\left(\frac{1}{2}(2L+3)\right)^{2}}z^{L-1}\,_{1}\tilde{F}_{2}\left(\frac{L}{2}+1;L+\frac{3}{2},\frac{L}{2}+\frac{3}{2};\frac{z}{4}\right)^{2}\end{array}\quad.\label{eq:final0_11}\end{equation}

Since the parameters in the $\,_{3}F_{4}$ hypergeometric function
contain $\frac{\mu}{2}+\frac{\nu}{2}$ in most of its parameters,
we can simplify the parameters on the left-hand side somewhat by letting $\mu\to a-\nu-1$ and $\nu\to b-1$
so that

\begin{equation}
\begin{array}{ccc}
 & &  \hspace{-8.4cm}   \,_{3}F_{4} \left(\frac{a}{2}+\frac{1}{2},\frac{a}{2},\frac{a}{2};\frac{a}{2}+1,a,a-b+1,b;z\right)   = a\Gamma(b)\Gamma(a-b+1)\\
 & \times & {\displaystyle \sum_{L=0}^{\infty}\frac{\pi^{2}(2L+1)i^{-2a+4L+2}2^{a-8L-7}\left((-1)^{b+L-1}+1\right)\Gamma(2L+2)^{2}\left((-1)^{a-b+L}+1\right)}{\Gamma\left(L+\frac{3}{2}\right)^{2}\Gamma\left(\frac{1}{2}(2L+3)\right)^{2}\Gamma\left(-\frac{b}{2}+\frac{L}{2}+\frac{3}{2}\right)\Gamma\left(\frac{b}{2}+\frac{L}{2}+\frac{1}{2}\right)\Gamma\left(\frac{a}{2}-\frac{b}{2}+\frac{L}{2}+1\right)\Gamma\left(-\frac{a}{2}+\frac{b}{2}+\frac{L}{2}+1\right)} }\\
 & \times & z^{\frac{1-a}{2}+L}   \,_{2}F_{3}\left(\frac{L}{2}+\frac{1}{2},\frac{L}{2}+1;-\frac{b}{2}+\frac{L}{2}+\frac{3}{2},\frac{b}{2}+\frac{L}{2}+\frac{1}{2},L+\frac{3}{2};\frac{z}{4}\right)  \\
 & \times &  \,_{2}F_{3}\left(\frac{L}{2}+\frac{1}{2},\frac{L}{2}+1;\frac{a}{2}-\frac{b}{2}+\frac{L}{2}+1,-\frac{a}{2}+\frac{b}{2}+\frac{L}{2}+1,L+\frac{3}{2};\frac{z}{4}\right)\\
 & = & a\Gamma(b)\Gamma(a-b+1)\\
 & \times & {\displaystyle \sum_{L=0}^{\infty}\frac{\pi^{2}(2L+1)i^{-2a+4L+2}2^{a-8L-7}\left((-1)^{b+L-1}+1\right)\Gamma(2L+2)^{2}\left((-1)^{a-b+L}+1\right)}{\Gamma\left(\frac{1}{2}(2L+3)\right)^{2}}} \\
 & \times & z^{\frac{1-a}{2}+L} \,_{2}\tilde{F}_{3}\left(\frac{L}{2}+\frac{1}{2},\frac{L}{2}+1;L+\frac{3}{2},-\frac{b}{2}+\frac{L}{2}+\frac{3}{2},\frac{b}{2}+\frac{L}{2}+\frac{1}{2};\frac{z}{4}\right)  \\
 & \times & \,_{2}\tilde{F}_{3}\left(\frac{L}{2}+\frac{1}{2},\frac{L}{2}+1;L+\frac{3}{2},-\frac{a}{2}+\frac{b}{2}+\frac{L}{2}+1,\frac{a}{2}-\frac{b}{2}+\frac{L}{2}+1;\frac{z}{4}\right)\end{array}\quad.\label{eq:final0_ab}\end{equation}
In this form, we must have $a=1,\,3,\,5,\,7\cdots$ and $ $$b=1,\,2,\,3,\cdots,\, a$
.

\section{For p=1}

Consider, now, an additional power of the cosine function (which equals $x$ in what follows)  that
will give a different Fourier-Legendre series than does a pure Bessel
function. The key step in modifying the prior derivation,  generalizes
the expansion coefficient at one point to

\begin{equation}
a_{LN}^{p}\left(k\right)=\frac{2L+1}{2}\frac{i^{-N}}{\pi}\int_{0}^{\pi}\left[\int_{-1}^{1}e^{ikx\cos\theta}x^{p}P_{L}(x)dx\right]\cos\left(N\theta\right)d\theta\quad,\label{eq:scriptjbeforelpinsertion}\end{equation}
in which we use \cite{GR5 p. 987 No. 8.511.4,Joachain p. 671 Eq. (B.46)}

\begin{equation}
e^{ikx\cos\theta}=\sum_{l'}(2l'+1)i^{l'}j_{l'}\left(k\cos\theta\right)P_{l'}(x)\quad.\label{eq:GR5 p. 987 No. 8.511.4}\end{equation}

\noindent Then using the recurrence relation for Legendre functions\cite{Joachain p. 666 eq. (B.4)}
one finds

\begin{equation}
\begin{array}{ccc}
a_{LN}^{p}\left(k\right) & = & \frac{2L+1}{2}\frac{i^{-N}}{\pi}\int_{0}^{\pi}\left[\int_{-1}^{1}\left(\sum_{l'=0}^{\infty}(2l'+1)i^{l'}j_{l'}\left(k\cos\theta\right)P_{l'}(x)\right)x^{p}P_{L}(x)dx\right]\cos\left(N\theta\right)d\theta\\
a_{LN}^{1}\left(k\right) & = & \frac{2L+1}{2}\frac{i^{-N}}{\pi}\int_{0}^{\pi}\left[\int_{-1}^{1}\left(\sum_{l'=0}^{\infty}(2l'+1)i^{l'}j_{l'}\left(k\cos\theta\right)\right)\left\{ \frac{l'P_{l'-1}(x)}{2l'+1}+\frac{\left(l'+1\right)P_{l'+1}(x)}{2l'+1}\right\} P_{L}(x)dx\right] \\
 & \times & \cos\left(N\theta\right)d\theta\quad,\end{array}\label{eq:coefp}\end{equation}
which complicates the use of Legendre function orthogonality\cite{Joachain p. 666 eq. (B.5)}
to truncate the series. For $l=0$ only the left-hand term in the
curly brackets will contribute with the $l'=1$ term in the sum. When
$l=1$, the left-hand term will again contribute with the $l'=2$
term in the sum, as will the right-hand term with the $l'=0$ term
in the sum. This pattern continues with the first term nonzero only
when $l'=l+1$, and the second nonzero only when $l'=l-1$. Thus

\begin{equation}
\begin{array}{ccc}
a_{LN}^{1}\left(k\right) & = & \frac{2L+1}{2}\frac{i^{-N}}{\pi}\int_{0}^{\pi}\left[\int_{-1}^{1}\left(\sum_{l'=0}^{\infty}i^{l'}j_{l'}\left(k\cos\theta\right)\right)\left\{ \frac{2l'+1}{2l'+1}\right\} \left\{ l'P_{l'-1}(x)+\left(l'+1\right)P_{l'+1}(x)\right\} P_{L}(x)dx\right]  \\
 & \times &\cos\left(N\theta\right)d\theta\\
 & = & \frac{2L+1}{2}\frac{i^{-N}}{\pi}\int_{0}^{\pi}\left[\int_{-1}^{1}\left(\sum_{l'=0}^{\infty}i^{l'}j_{l'}\left(k\cos\theta\right)\right)\left\{ \frac{2l'}{2l'-1}\delta_{l'-1,L}+\frac{2\left(l'+1\right)}{2l'+3}\delta_{l'+1,L}\right\} \right]\cos\left(N\theta\right)d\theta\\
 & = & \frac{2L+1}{2}\frac{i^{-N}}{\pi}\int_{0}^{\pi}\left[{\displaystyle \sum_{\ell=0}^{\infty}i^{\ell+1}j_{\ell+1}\left(k\cos\theta\right)\frac{2\left(\ell+1\right)}{2\ell+1}\delta_{\ell L}+\sum_{\ell=1}^{\infty}i^{\ell-1}j_{\ell-1}\left(k\cos\theta\right)\frac{2\ell}{2\ell+1}\delta_{\ell L}}\right]  \\
 & \times & \cos\left(N\theta\right)d\theta\\
 & = & \frac{i^{-N}}{\pi}\int_{0}^{\pi}\left[i^{L+1}j_{L+1}\left(k\cos\theta\right)\left(L+1\right)+i^{L-1}j_{L-1}\left(k\cos\theta\right)L\right]\cos\left(N\theta\right)d\theta\end{array}\label{eq:coefp=00003D1}
 \end{equation}
 where the factor of \emph{L} takes the second term to zero when $L=0$
so we need not worry about spherical Bessel functions with an index
less than zero. (This is the same reason why the first term of the
$l'$ sum in the second line really starts at $1$ so that the $\ell$
sum starts at $0$ for the first sum in the third line.) Note that
\cite{GH p. 23 No. 171.4} gives only the second of these two terms,
but one may interchange the meaning of \emph{m} and \emph{n} in this
to give the first, the central line in\begin{equation}
\int_{-1}^{1}xP_{n}(x)P_{m}(x)\, dx=\begin{cases}
\frac{2n+2}{(2n+1)(2n+3)} & m=n+1\\
\frac{2n}{(2n-1)(2n+1)} & m=n-1\\
0 & otherwise\end{cases} \quad . \label{eq:GH p. 23 No. 171.4extended}\end{equation}

Using the series expansion \cite{wolfram.com/03.21.06.0019.01}

\begin{equation}
j_{L\pm1}(x)=\sqrt{\pi}\,2^{-(L\pm1)-1}x^{L\pm1}\sum_{M=0}^{\infty}\frac{\left(-\frac{1}{4}\right)^{M}x^{2M}}{M!\Gamma\left(M+(L\pm1)+\frac{3}{2}\right)}\label{eq:spherbesselj}\end{equation}
 we find

\begin{equation}
\begin{array}{ccc}
a_{LN}^{1} & \left(k\right)= & \frac{i^{-N}}{\pi}\int_{0}^{\pi}\left[i^{L+1}j_{L+1}\left(k\alpha_{0}\cos\theta\right)\left(L+1\right)+i^{L-1}j_{L-1}\left(k\alpha_{0}\cos\theta\right)L\right]\cos\left(N\theta\right)d\theta\\
 & = & \frac{1}{2\sqrt{\pi}}i^{-N}\int_{0}^{\pi}\left[\left(\frac{i}{2}\right)^{L+1}\left(L+1\right){\displaystyle \sum_{M=0}^{\infty}\frac{\left(-1\right)^{M}}{4^{M}M!\Gamma\left(L+M+\frac{5}{2}\right)}}\left(k\cos(\theta)\right){}^{2M+L+1}\right.\\
 & + & \left.\left(\frac{i}{2}\right)^{L-1}L{\displaystyle \sum_{M=0}^{\infty}\frac{\left(-1\right)^{M}}{4^{M}M!\Gamma\left(L+M+\frac{1}{2}\right)}}\left(k\cos(\theta)\right){}^{2M+L-1}\right]\cos\left(N\theta\right)d\theta\end{array}\label{eq:scriptj1Msum}\end{equation}

This integral may be done using an integral -- that has three branches
over the interval $\left[0,\frac{\pi}{2}\right]$ -- that Gr\"{o}bner
and Hofreiter\cite{GH p. 110 No. 332.14a} extended to the interval
$\left[0,\pi\right]$ with a prefactor $(1+(-1)^{m+n})$, which renders
the central one of the three possibilities nonzero only for even values
for $m+n$.

\begin{equation}
\int_{0}^{\pi}cos^{m}\theta\, cos(n\theta)d\theta=\left(1+(-1)^{m+n}\right)\frac{\pi}{2^{m+1}}\binom{m}{\frac{m-n}{2}}\quad\left[m\geq n>-1,\, m-n=2K\right]\quad.\label{eq:GR5 p. 417 No. 3.631.17}\end{equation}
 The other two branches, for odd $m+n$ on $\left[0,\frac{\pi}{2}\right]$,
are zero on $\left[0,\pi\right]$ when this prefactor is included.
(The lower limit on \emph{m} $\left[m\geq n>-1\right]$ was a finding
of the prior paper.) Then

\begin{equation}
\begin{array}{ccc}
a_{LN}^{1}\left(k\right) & = & \frac{1}{2\sqrt{\pi}}i^{-N}\left[\left(\frac{i}{2}\right)^{L+1}\left(L+1\right){\displaystyle \sum_{M=0}^{\infty}\frac{\left(-1\right)^{M}}{4^{M}M!\Gamma\left(L+M+\frac{5}{2}\right)}}\left(k\right){}^{2M+L+1}\left(1+(-1)^{2M+L+1+N}\right) \right.\\
 & \times & \frac{\pi}{2^{2M+L+2}}\binom{2M+L+1}{\frac{1}{2}(2M+L+1-N)}\\
 & + & \left.\left(\frac{i}{2}\right)^{L-1}L{\displaystyle \sum_{M=0}^{\infty}\frac{\left(-1\right)^{M}}{4^{M}M!\Gamma\left(L+M+\frac{1}{2}\right)}}\left(k\right){}^{2M+L-1}\left(1+(-1)^{2M+L-1+N}\right)  \right.\\
 & \times & \left.  \frac{\pi}{2^{2M+L}}  \binom{2M+L-1}{\frac{1}{2}(2M+L-1-N)}\right]\\
 & = & \hspace{-4.9cm}  \frac{\sqrt{\pi}}{2}i^{-N}\left[\frac{i^{L-1}2^{1-2L}L\, k^{L-1}}{\Gamma\left(\frac{1}{2}(2L+1)\right)}\left((-1)^{L+N+1}+1\right)\binom{L-1}{\frac{1}{2}(L-N-1)} \right. \\
 & \times & \left. \,_{2}F_{3}\left(\frac{L}{2}+\frac{1}{2},\frac{L}{2};L+\frac{1}{2},\frac{L}{2}-\frac{N}{2}+\frac{1}{2},\frac{L}{2}+\frac{N}{2}+\frac{1}{2};-\frac{k^{2}}{4}\right)\right.\\
 & + &  \hspace{-2.8cm} \left.\frac{i^{L+1}2^{-2L-3}(L+1)k^{L+1}}{\Gamma\left(\frac{1}{2}(2L+5)\right)}\left(1+(-1)^{L+N-1}\right)\binom{L+1}{\frac{1}{2}(L-N+1)}  \right. \\
 & \times & \left.  \,_{2}F_{3}\left(\frac{L}{2}+1,\frac{L}{2}+\frac{3}{2};L+\frac{5}{2},\frac{L}{2}-\frac{N}{2}+\frac{3}{2},\frac{L}{2}+\frac{N}{2}+\frac{3}{2};-\frac{k^{2}}{4}\right)\right]\\
 & = &  \hspace{-5.4cm} \frac{\sqrt{\pi}}{2}i^{-N}\left[i^{L-1}2^{1-2L}L\, k^{L-1}\left(1+(-1)^{L+N+1}\right)\Gamma(L)   \right. \\
 & \times & \left.  \,_{2}\tilde{F}_{3}\left(\frac{L}{2}+\frac{1}{2},\frac{L}{2};L+\frac{1}{2},\frac{L}{2}-\frac{N}{2}+\frac{1}{2},\frac{L}{2}+\frac{N}{2}+\frac{1}{2};-\frac{k^{2}}{4}\right)\right.\\
 & + &  \hspace{-2.5cm} i^{L+1}2^{-2L-3}(L+1)k^{L+1}\left(1+(-1)^{L+N-1}\right)\Gamma(L+2)\\
& \times &  \left. \,_{2}\tilde{F}_{3}\left(\frac{L}{2}+1,\frac{L}{2}+\frac{3}{2};L+\frac{5}{2},\frac{L}{2}-\frac{N}{2}+\frac{3}{2},\frac{L}{2}+\frac{N}{2}+\frac{3}{2};-\frac{k^{2}}{4}\right)\right]\end{array}\label{eq:a1LN_p. 110 No. 332.14a}\end{equation}

For the special cases $N=\{0,1\}$these reduce somewhat \cite{wolfram.com/07.26.03.0002.01}
to

\begin{equation}
\begin{array}{ccc}
a_{L0}^{1}\left(k\right) & = & \frac{\sqrt{\pi}}{2}\left[\frac{i^{L+1}2^{-2L-3}}{\Gamma\left(\frac{L}{2}+\frac{3}{2}\right)}\left(1+(-1)^{L-1}\right)(L+1)k^{L+1}\Gamma(L+2)\,_{1}\tilde{F}_{2}\left(\frac{L}{2}+1;L+\frac{5}{2},\frac{L}{2}+\frac{3}{2};-\frac{k^{2}}{4}\right)\right.\\
 & + & \left.\frac{i^{L-1}2^{1-2L}}{\Gamma\left(\frac{L}{2}+\frac{1}{2}\right)}\left(1+(-1)^{L+1}\right)Lk^{L-1}\Gamma(L)\,_{1}\tilde{F}_{2}\left(\frac{L}{2};L+\frac{1}{2},\frac{L}{2}+\frac{1}{2};-\frac{k^{2}}{4}\right)\right]\\
a_{L1}^{1}\left(k\right) & = & \frac{\sqrt{\pi}}{2}(-i)\left[\frac{i^{L-1}2^{1-2L}}{\Gamma\left(\frac{L}{2}\right)}\left(1+(-1)^{L+2}\right)Lk^{L-1}\Gamma(L)\,_{1}\tilde{F}_{2}\left(\frac{L}{2}+\frac{1}{2};L+\frac{1}{2},\frac{L}{2}+1;-\frac{k^{2}}{4}\right)\right.\quad.\\
 & + & \left.\frac{i^{L+1}2^{-2L-3}}{\Gamma\left(\frac{L}{2}+1\right)}\left(1+(-1)^{L}\right)(L+1)k^{L+1}\Gamma(L+2)\,_{1}\tilde{F}_{2}\left(\frac{L}{2}+\frac{3}{2};L+\frac{5}{2},\frac{L}{2}+2;-\frac{k^{2}}{4}\right)\right]\end{array}\label{eq:a1L0_1p=00003D1}\end{equation}

We will not display intermediate results as we step outward since
they mirror those for $p=0$ but with twice as many terms. With $x=k\alpha_{0}$,
the final form is

\begin{equation}
\begin{array}{ccc}
\Xi_{n}^{1} & & \hspace{-9.4cm}  \left(k\alpha_{0},-\frac{z}{2}\right)  = \sum_{l0}F_{l0}^{1n}(k\alpha_{0},-\frac{z}{2})F_{l0}^{1n}(k\alpha_{0},-\frac{z}{2})\\
 & = & \sum_{l0}\frac{2l+1}{4\pi}\left[{\displaystyle \sum_{j=-\left(n-\delta\right)/2}^{\infty}\quad}{\displaystyle \sum_{M=-\left(n-\delta\right)/2}^{\infty}}J_{-j}\left(-\frac{z}{2}\right)J_{-M}\left(-\frac{z}{2}\right)\right.\\
 & \times & {\displaystyle \frac{\pi^{3}2^{4-4l}x^{2l-2}}{(2l+1)^{2}}}\left((-1)^{2j+l+n+1}+1\right)\left((-1)^{l+2M+n+1}+1\right)(-1)^{-j+l-M-n-1}\\
 & \times & \left(\frac{l}{\Gamma\left(\frac{1}{2}(2l+1)\right)}\binom{l-1}{\frac{1}{2}(-2j+l-n-1)}\,_{2}F_{3}\left(\frac{l}{2}+\frac{1}{2},\frac{l}{2};l+\frac{1}{2},-j+\frac{l}{2}-\frac{n}{2}+\frac{1}{2},j+\frac{l}{2}+\frac{n}{2}+\frac{1}{2};-\frac{x^{2}}{4}\right)\right.\\
 & - & \left.\frac{(l+1)x^{2}}{16\Gamma\left(\frac{1}{2}(2l+5)\right)}\binom{l+1}{\frac{1}{2}(-2j+l-n+1)}\,_{2}F_{3}\left(\frac{l}{2}+1,\frac{l}{2}+\frac{3}{2};l+\frac{5}{2},-j+\frac{l}{2}-\frac{n}{2}+\frac{3}{2},j+\frac{l}{2}+\frac{n}{2}+\frac{3}{2};-\frac{x^{2}}{4}\right)\right)\\
 & \times & \left(\frac{l}{\Gamma\left(\frac{1}{2}(2l+1)\right)}\binom{l-1}{\frac{1}{2}(l-2M-n-1)}\,_{2}F_{3}\left(\frac{l}{2}+\frac{1}{2},\frac{l}{2};l+\frac{1}{2},\frac{l}{2}-M-\frac{n}{2}+\frac{1}{2},\frac{l}{2}+M+\frac{n}{2}+\frac{1}{2};-\frac{x^{2}}{4}\right)\right.\\
 & - & \left.\frac{(l+1)x^{2}}{16\Gamma\left(\frac{1}{2}(2l+5)\right)}\binom{l+1}{\frac{1}{2}(l-2M-n+1)}\,_{2}F_{3}\left(\frac{l}{2}+1,\frac{l}{2}+\frac{3}{2};l+\frac{5}{2},\frac{l}{2}-M-\frac{n}{2}+\frac{3}{2},\frac{l}{2}+M+\frac{n}{2}+\frac{3}{2};-\frac{x^{2}}{4}\right)\right)\\
 & + & 2\left(-1\right)^{n}{\displaystyle \sum_{j=\left(n-\delta\right)/2+1}^{\infty}}{\displaystyle \sum_{M=-\left(n-\delta\right)/2}^{\infty}}J_{j}\left(-\frac{z}{2}\right)J_{-M}\left(-\frac{z}{2}\right)\\
 & \times & \frac{\pi^{3}2^{4-4l}x^{2l-2}}{(2l+1)^{2}}(-1)^{-j+l-M-1}\left((-1)^{2j+l-n+1}+1\right)\left((-1)^{l+2M+n+1}+1\right)\\
 & \times & \left(\frac{l}{\Gamma\left(\frac{1}{2}(2l+1)\right)}\binom{l-1}{\frac{1}{2}(-2j+l+n-1)}\,_{2}F_{3}\left(\frac{l}{2}+\frac{1}{2},\frac{l}{2};l+\frac{1}{2},j+\frac{l}{2}-\frac{n}{2}+\frac{1}{2},-j+\frac{l}{2}+\frac{n}{2}+\frac{1}{2};-\frac{x^{2}}{4}\right)\right.\\
 & - & \left.\frac{(l+1)x^{2}}{16\Gamma\left(\frac{1}{2}(2l+5)\right)}\binom{l+1}{\frac{1}{2}(-2j+l+n+1)}\,_{2}F_{3}\left(\frac{l}{2}+1,\frac{l}{2}+\frac{3}{2};l+\frac{5}{2},j+\frac{l}{2}-\frac{n}{2}+\frac{3}{2},-j+\frac{l}{2}+\frac{n}{2}+\frac{3}{2};-\frac{x^{2}}{4}\right)\right)\\
 & \times & \left(\frac{l}{\Gamma\left(\frac{1}{2}(2l+1)\right)}\binom{l-1}{\frac{1}{2}(l-2M-n-1)}\,_{2}F_{3}\left(\frac{l}{2}+\frac{1}{2},\frac{l}{2};l+\frac{1}{2},\frac{l}{2}-M-\frac{n}{2}+\frac{1}{2},\frac{l}{2}+M+\frac{n}{2}+\frac{1}{2};-\frac{x^{2}}{4}\right)\right.\\
 & - & \left.\frac{(l+1)x^{2}}{16\Gamma\left(\frac{1}{2}(2l+5)\right)}\binom{l+1}{\frac{1}{2}(l-2M-n+1)}\,_{2}F_{3}\left(\frac{l}{2}+1,\frac{l}{2}+\frac{3}{2};l+\frac{5}{2},\frac{l}{2}-M-\frac{n}{2}+\frac{3}{2},\frac{l}{2}+M+\frac{n}{2}+\frac{3}{2};-\frac{x^{2}}{4}\right)\right)\\
 & + & {\displaystyle \sum_{j=\left(n-\delta\right)/2+1}^{\infty}}{\displaystyle {\displaystyle \sum_{M=\left(n-\delta\right)/2+1}^{\infty}}}J_{j}\left(-\frac{z}{2}\right)J_{M}\left(-\frac{z}{2}\right)\\
 & \times & {\displaystyle \frac{\pi^{3}2^{4-4l}x^{2l-2}}{(2l+1)^{2}}}\left((-1)^{2j+l-n+1}+1\right)\left((-1)^{l+2M-n+1}+1\right)(-1)^{-j+l-M+n-1}\\
 & \times & \left(\frac{l}{\Gamma\left(\frac{1}{2}(2l+1)\right)}\binom{l-1}{\frac{1}{2}(-2j+l+n-1)}\,_{2}F_{3}\left(\frac{l}{2}+\frac{1}{2},\frac{l}{2};l+\frac{1}{2},j+\frac{l}{2}-\frac{n}{2}+\frac{1}{2},-j+\frac{l}{2}+\frac{n}{2}+\frac{1}{2};-\frac{x^{2}}{4}\right)\right.\\
 & - & \left.\frac{(l+1)x^{2}}{16\Gamma\left(\frac{1}{2}(2l+5)\right)}\binom{l+1}{\frac{1}{2}(-2j+l+n+1)}\,_{2}F_{3}\left(\frac{l}{2}+1,\frac{l}{2}+\frac{3}{2};l+\frac{5}{2},j+\frac{l}{2}-\frac{n}{2}+\frac{3}{2},-j+\frac{l}{2}+\frac{n}{2}+\frac{3}{2};-\frac{x^{2}}{4}\right)\right)\\
 & \times & \left(\frac{l}{\Gamma\left(\frac{1}{2}(2l+1)\right)}\binom{l-1}{\frac{1}{2}(l-2M+n-1)}\,_{2}F_{3}\left(\frac{l}{2}+\frac{1}{2},\frac{l}{2};l+\frac{1}{2},\frac{l}{2}+M-\frac{n}{2}+\frac{1}{2},\frac{l}{2}-M+\frac{n}{2}+\frac{1}{2};-\frac{x^{2}}{4}\right)\right.\\
 & - & \left.\left.\frac{(l+1)x^{2}}{16\Gamma\left(\frac{1}{2}(2l+5)\right)}\binom{l+1}{\frac{1}{2}(l-2M+n+1)}\,_{2}F_{3}\left(\frac{l}{2}+1,\frac{l}{2}+\frac{3}{2};l+\frac{5}{2},\frac{l}{2}+M-\frac{n}{2}+\frac{3}{2},\frac{l}{2}-M+\frac{n}{2}+\frac{3}{2};-\frac{x^{2}}{4}\right)\right)\right]\end{array}\label{eq:chi0sump=00003D00003D1}\end{equation}

On the other hand, for $p=1$ in the direct integration method, the
only thing that changes is the power of the $\cos^{2}\left(\theta\right)=u^{2}=y$
factor in

\begin{eqnarray}
\Xi_{n}^{1}\left(k\alpha_{0},-\frac{z}{2}\right) & = & \int\cos^{2}\left(\theta\right)J_{n}^{2}\left(\mathbf{k}\cdot\boldsymbol{\alpha}_{0},-\frac{z}{2}\right)d\Omega=2\pi\int_{-\pi}^{\pi}\cos^{2}\left(\theta\right)J_{n}^{2}\left(k\alpha_{0}\cos(\theta),-\frac{z}{2}\right)d\left(\cos(\theta)\right) \nonumber \\
& = & 2\pi\int_{-1}^{1}f\left(u^{2}\right)du=4\pi\int_{0}^{1}f\left(u^{2}\right)du=2\pi\int_{0}^{1}f\left(y\right)y^{-1/2}dy\quad.\label{eq:angintsetup-1}\end{eqnarray}
 and in $\alpha=\left\{ j+M+n+1,\: j+M+1,\: j+M+1,\: j+M-n+1\right\} +p$,
in the following 

\begin{eqnarray}
\int_{0}^{a}& & \hspace{-1.0cm} y^{\alpha-1}  (a-y)^{\beta-1}\,_{p}F_{q}(a_{1},\ldots,\, a_{p},b_{1},\ldots,\, b_{q};-\omega y)\, dy \\
& = & \frac{\Gamma(\alpha)\Gamma(\beta)a^{\alpha+\beta-1}}{\Gamma(\alpha+\beta)}\,_{p+1}F_{q+1}(a_{1},\ldots,\, a_{p},\alpha;b_{1},\ldots,\, b_{q},\alpha+\beta;-a\omega)\quad\left[\Re(\alpha)>0\land\Re(\beta)>0\land a>0\right] \nonumber \label{eq:PBM3 p. 334 No. 2.22.2.1-1}\end{eqnarray}
 %
 that simply raises $\alpha$ by one for each of the four terms; with
$a=1$ and $\beta=1$ unchanged.

Then the above and any of the four terms in equation (\ref{eq:chi0sump=00003D00003D1})
give

\begin{equation}
\begin{array}{ccc}
\,_{3}F_{4} & & \hspace{-1.3cm}  \left(\frac{\mu}{2}+\frac{\nu}{2}+\frac{1}{2},\frac{\mu}{2}+\frac{\nu}{2}+1,\frac{\mu}{2}+\frac{\nu}{2}+\frac{1}{2}+p;\mu+1,\nu+1,\mu+\nu+1,\frac{\mu}{2}+\frac{\nu}{2}+\frac{3}{2}+p;z\right)_{p\rightarrow1} \\
& = & \pi i^{-\mu-\nu}2^{\mu+\nu-2}\Gamma(\mu+1)\Gamma(\nu+1)(2p+\mu+\nu+1)z^{\frac{1}{2}(-\mu-\nu)}\\
 & \times & {\displaystyle \sum_{L=0}^{\infty}\frac{i^{2L-2}2^{2-4L}\left(1+(-1)^{L+\mu+1}\right)\left(1+(-1)^{L+\nu+1}\right)(-1)^{L-\frac{\mu}{2}-\frac{\nu}{2}-1}}{2L+1}}z^{L-1}\\
 & \times & \left(\frac{L}{\Gamma\left(\frac{1}{2}(2L+1)\right)}\binom{L-1}{\frac{1}{2}(L-\mu-1)}\,_{2}F_{3}\left(\frac{L}{2}+\frac{1}{2},\frac{L}{2};L+\frac{1}{2},\frac{L}{2}-\frac{\mu}{2}+\frac{1}{2},\frac{L}{2}+\frac{\mu}{2}+\frac{1}{2};\frac{z}{4}\right)\right.\\
 & + & \left.\frac{(L+1)}{16\Gamma\left(\frac{1}{2}(2L+5)\right)}z\binom{L+1}{\frac{1}{2}(L-\mu+1)}\,_{2}F_{3}\left(\frac{L}{2}+1,\frac{L}{2}+\frac{3}{2};L+\frac{5}{2},\frac{L}{2}-\frac{\mu}{2}+\frac{3}{2},\frac{L}{2}+\frac{\mu}{2}+\frac{3}{2};\frac{z}{4}\right)\right)\\
 & \times & \left(\frac{L}{\Gamma\left(\frac{1}{2}(2L+1)\right)}\binom{L-1}{\frac{1}{2}(L-\nu-1)}\,_{2}F_{3}\left(\frac{L}{2}+\frac{1}{2},\frac{L}{2};L+\frac{1}{2},\frac{L}{2}-\frac{\nu}{2}+\frac{1}{2},\frac{L}{2}+\frac{\nu}{2}+\frac{1}{2};\frac{z}{4}\right)\right.\\
 & + & \left.\frac{(L+1)}{16\Gamma\left(\frac{1}{2}(2L+5)\right)}z\binom{L+1}{\frac{1}{2}(L-\nu+1)}\,_{2}F_{3}\left(\frac{L}{2}+1,\frac{L}{2}+\frac{3}{2};L+\frac{5}{2},\frac{L}{2}-\frac{\nu}{2}+\frac{3}{2},\frac{L}{2}+\frac{\nu}{2}+\frac{3}{2};\frac{z}{4}\right)\right)\\
 & = & \pi i^{-\mu-\nu}2^{\mu+\nu-2}\Gamma(\mu+1)\Gamma(\nu+1)(2p+\mu+\nu+1)z^{\frac{1}{2}(-\mu-\nu)}\\
 & \times & {\displaystyle \sum_{L=0}^{\infty}\frac{i^{2L-2}2^{2-4L}\left(1+(-1)^{L+\mu+1}\right)\left(1+(-1)^{L+\nu+1}\right)(-1)^{L-\frac{\mu}{2}-\frac{\nu}{2}-1}}{2L+1}}z^{L-1}\\
 & \times & \left(\Gamma(L+1)\,_{2}\tilde{F}_{3}\left(\frac{L}{2}+\frac{1}{2},\frac{L}{2};L+\frac{1}{2},\frac{L}{2}-\frac{\mu}{2}+\frac{1}{2},\frac{L}{2}+\frac{\mu}{2}+\frac{1}{2};\frac{z}{4}\right)\right.\\
 & + & \left.\frac{1}{16}(L+1)\Gamma(L+2)\, z\,_{2}\tilde{F}_{3}\left(\frac{L}{2}+1,\frac{L}{2}+\frac{3}{2};L+\frac{5}{2},\frac{L}{2}-\frac{\mu}{2}+\frac{3}{2},\frac{L}{2}+\frac{\mu}{2}+\frac{3}{2};\frac{z}{4}\right)\right)\\
 & \times & \left(\Gamma(L+1)\,_{2}\tilde{F}_{3}\left(\frac{L}{2}+\frac{1}{2},\frac{L}{2};L+\frac{1}{2},\frac{L}{2}-\frac{\nu}{2}+\frac{1}{2},\frac{L}{2}+\frac{\nu}{2}+\frac{1}{2};\frac{z}{4}\right)\right.\\
 & + & \left.\frac{1}{16}(L+1)\Gamma(L+2)\, z\,_{2}\tilde{F}_{3}\left(\frac{L}{2}+1,\frac{L}{2}+\frac{3}{2};L+\frac{5}{2},\frac{L}{2}-\frac{\nu}{2}+\frac{3}{2},\frac{L}{2}+\frac{\nu}{2}+\frac{3}{2};\frac{z}{4}\right)\right)\end{array}\quad.\label{eq:final0-1}\end{equation}

For $p=1$ the right-hand side requires  four nonzero
terms in the sum to get seven-digit accuracy for $\mu+\nu<10$,

\begin{tabular}{|l|l|c|c|c|}
\hline 
Left-hand side & Right-hand side & $\mu$ & $\nu$ & z\tabularnewline
\hline
\hline 
1.052175485266236 & 1.052175485266234 & 0 & 0 & 0.17\tabularnewline
\hline 
1.01448216 & 1.01448208 & 4 & 4 & 0.17\tabularnewline
\hline 
1.0307786670736816 & 1.03077866707368164 & 1 & 1 & 0.17\tabularnewline
\hline 
1.01234929 & 1.01234927 & 5 & 5 & 0.17\tabularnewline
\hline 
\end{tabular}

Table 2. The left and right sides of (\ref{eq:final0-1}) when the
variable is set to the arbitrary value $z=0.17$ and we include four
nonzero terms in the sum, shown through the digit with which the two
sides disagree. 

\vspace{0.3cm}

When $\mu=\nu=0$ the order of the hypergeometric functions is reduced
since the parameters $a_{2}=b_{1}$ on the left-hand side. On the
right-hand side, $a_{1}=b_{2}$ in the first hypergeometric function
in the summed pair within the products and $a_{2}=b_{2}$ in the second,
giving the special case,

\begin{equation}
\begin{array}{ccc}
\,_{2}F_{3} & & \hspace{-6.3cm}  \left(\frac{1}{2},p+\frac{1}{2};1,1,p+\frac{3}{2};z\right)_{p\rightarrow1} \\
& = & \frac{1}{4}\pi(2p+1)_{p\rightarrow1}\sum_{L=0}^{\infty}\frac{(-1)^{L-1}i^{2L-2}2^{2-4L}\left(1+(-1)^{L+1}\right)^{2}}{2L+1}z^{L-1}  \\
& \times &  \left(\frac{(L+1)}{16\Gamma\left(\frac{1}{2}(2L+5)\right)}\binom{L+1}{\frac{L+1}{2}}\, z\,_{1}F_{2}\left(\frac{L}{2}+1;\frac{L}{2}+\frac{3}{2},L+\frac{5}{2};\frac{z}{4}\right)\right.\\
 & + & \left|\frac{L}{\Gamma\left(\frac{1}{2}(2L+1)\right)}\binom{L-1}{\frac{L-1}{2}}\,_{1}F_{2}\left(\frac{L}{2};\frac{L}{2}+\frac{1}{2},L+\frac{1}{2};\frac{z}{4}\right)\right)^{2}\\
 & = & \frac{1}{4}\pi(2p+1)_{p\rightarrow1}\sum_{L=0}^{\infty}\frac{(-1)^{L-1}i^{2L-2}2^{2-4L}\left(1+(-1)^{L+1}\right)^{2}}{2L+1}z^{L-1}   \\
& \times & \left(\frac{(L+1)\Gamma(L+2)}{16\Gamma\left(\frac{L}{2}+\frac{3}{2}\right)}\, z\,_{1}\tilde{F}_{2}\left(\frac{L}{2}+1;L+\frac{5}{2},\frac{L}{2}+\frac{3}{2};\frac{z}{4}\right)\right.\\
 & + & \left|\frac{\Gamma(L+1)}{\Gamma\left(\frac{L}{2}+\frac{1}{2}\right)}\,_{1}\tilde{F}_{2}\left(\frac{L}{2};L+\frac{1}{2},\frac{L}{2}+\frac{1}{2};\frac{z}{4}\right)\right)^{2}\end{array}\quad.\label{eq:final0_00-1}
 \end{equation}
 
When $\mu=\nu=1$ the order of the hypergeometric functions is also
reduced since the parameters $a_{2}=b_{1}$ on the left-hand side.
On the right-hand side, $a_{2}=b_{2}$ in the first hypergeometric
function in the summed pair within the products and $a_{1}=b_{2}$
in the second, giving the special case,

\begin{equation}
\begin{array}{ccc}
\,_{2}F_{3}  & & \hspace{-6.3cm}  \left(\frac{3}{2},p+\frac{3}{2};2,3,p+\frac{5}{2};z\right)_{p\rightarrow1} \\
& = & -\frac{\pi}{z}(2p+3)_{p\rightarrow1}\sum_{L=0}^{\infty}\frac{(-1)^{L-2}i^{2L-2}2^{2-4L}\left(1+(-1)^{L+2}\right)^{2}}{2L+1}z^{L-1} \\
& \times & \left(\frac{L}{\Gamma\left(\frac{1}{2}(2L+1)\right)}\binom{L-1}{\frac{L-2}{2}}\,_{1}F_{2}\left(\frac{L}{2}+\frac{1}{2};\frac{L}{2}+1,L+\frac{1}{2};\frac{z}{4}\right)\right.\\
 & + & \left|\frac{(L+1)}{16\Gamma\left(\frac{1}{2}(2L+5)\right)}z\binom{L+1}{\frac{L}{2}}\,_{1}F_{2}\left(\frac{L}{2}+\frac{3}{2};\frac{L}{2}+2,L+\frac{5}{2};\frac{z}{4}\right)\right)^{2}\\
 & = & -\frac{\pi}{z}(2p+3)_{p\rightarrow1}\sum_{L=0}^{\infty}\frac{(-1)^{L-2}i^{2L-2}2^{2-4L}\left(1(-1)^{L+2}\right)^{2}}{2L+1}z^{L-1} \\
& \times & \left(\frac{\Gamma(L+1)}{\Gamma\left(\frac{L}{2}\right)}\,_{1}\tilde{F}_{2}\left(\frac{L}{2}+\frac{1}{2};L+\frac{1}{2},\frac{L}{2}+1;\frac{z}{4}\right)\right.\\
 & + & \left|\frac{(L+1)\Gamma(L+2)}{16\Gamma\left(\frac{L}{2}+1\right)}\, z\,_{1}\tilde{F}_{2}\left(\frac{L}{2}+\frac{3}{2};L+\frac{5}{2},\frac{L}{2}+2;\frac{z}{4}\right)\right)^{2}\end{array}\quad.\label{eq:final0_11-1}\end{equation}

As with $p=0$, we can simplify the left-hand side of (\ref{eq:final0-1}) somewhat by letting
$\mu\to a-\nu-1$ and $\nu\to b-1$ so that$ $

\begin{equation}
\begin{array}{ccc}
\,_{3}F_{4}   & & \hspace{-8.8cm}  \left(\frac{a}{2}+\frac{1}{2},\frac{a}{2},\frac{a}{2}+p;a,a-b+1,b,\frac{a}{2}+p+1;z\right)_{p\rightarrow1} \\
& = & \pi i^{1-a}2^{a-3}z^{\frac{1-a}{2}}\Gamma(b)\Gamma(a-b+1)(a+2p)_{p\rightarrow1}\\
 & \times & {\displaystyle \sum_{L=0}^{\infty}\frac{i^{2L-2}2^{2-4L}\left(1+(-1)^{b+L}\right)(-1)^{\frac{b-a}{2}+\frac{1-b}{2}+L-1}\left(1+(-1)^{a-b+L+1}\right)}{2L+1}z^{L-1}}\\
 & \times & \left(\frac{L}{\Gamma\left(\frac{1}{2}(2L+1)\right)}\binom{L-1}{\frac{L-b}{2}}\,_{2}F_{3}\left(\frac{L}{2}+\frac{1}{2},\frac{L}{2};-\frac{b}{2}+\frac{L}{2}+1,\frac{b}{2}+\frac{L}{2},L+\frac{1}{2};\frac{z}{4}\right)\right.\\
 & + & \left.\frac{(L+1)}{16\Gamma\left(\frac{1}{2}(2L+5)\right)}\binom{L+1}{\frac{1}{2}(-b+L+2)}\, z\,_{2}F_{3}\left(\frac{L}{2}+1,\frac{L}{2}+\frac{3}{2};-\frac{b}{2}+\frac{L}{2}+2,\frac{b}{2}+\frac{L}{2}+1,L+\frac{5}{2};\frac{z}{4}\right)\right)\\
 & \times & \left(\frac{L}{\Gamma\left(\frac{1}{2}(2L+1)\right)}\binom{L-1}{\frac{1}{2}(-a+b+L-1)}\,_{2}F_{3}\left(\frac{L}{2}+\frac{1}{2},\frac{L}{2};\frac{a}{2}-\frac{b}{2}+\frac{L}{2}+\frac{1}{2},-\frac{a}{2}+\frac{b}{2}+\frac{L}{2}+\frac{1}{2},L+\frac{1}{2};\frac{z}{4}\right)\right.\\
 & + & \left.\frac{(L+1)}{16\Gamma\left(\frac{1}{2}(2L+5)\right)}\binom{L+1}{\frac{1}{2}(-a+b+L+1)}\, z\,_{2}F_{3}\left(\frac{L}{2}+1,\frac{L}{2}+\frac{3}{2};\frac{a}{2}-\frac{b}{2}+\frac{L}{2}+\frac{3}{2},-\frac{a}{2}+\frac{b}{2}+\frac{L}{2}+\frac{3}{2},L+\frac{5}{2};\frac{z}{4}\right)\right)\\
 & = & \pi i^{1-a}2^{a-3}z^{\frac{1-a}{2}}\Gamma(b)\Gamma(a-b+1)(a+2p)_{p\rightarrow1}\\
 & \times & {\displaystyle \sum_{L=0}^{\infty}\frac{i^{2L-2}2^{2-4L}\left(1+(-1)^{b+L}\right)(-1)^{\frac{b-a}{2}+\frac{1-b}{2}+L-1}\left(1+(-1)^{a-b+L+1}\right)}{2L+1}z^{L-1}}\\
 & \times & \left(\Gamma(L+1)\,_{2}\tilde{F}_{3}\left(\frac{L}{2}+\frac{1}{2},\frac{L}{2};L+\frac{1}{2},\frac{1-b}{2}+\frac{L}{2}+\frac{1}{2},\frac{b-1}{2}+\frac{L}{2}+\frac{1}{2};\frac{z}{4}\right)\right.\\
 & + & \left.\frac{1}{16}(L+1)z\Gamma(L+2)\,_{2}\tilde{F}_{3}\left(\frac{L}{2}+1,\frac{L}{2}+\frac{3}{2};L+\frac{5}{2},\frac{1-b}{2}+\frac{L}{2}+\frac{3}{2},\frac{b-1}{2}+\frac{L}{2}+\frac{3}{2};\frac{z}{4}\right)\right)\\
 & \times & \left(\Gamma(L+1)\,_{2}\tilde{F}_{3}\left(\frac{L}{2}+\frac{1}{2},\frac{L}{2};L+\frac{1}{2},\frac{b-a}{2}+\frac{L}{2}+\frac{1}{2},\frac{a-b}{2}+\frac{L}{2}+\frac{1}{2};\frac{z}{4}\right)\right.\\
 & + & \left.\frac{1}{16}(L+1)z\Gamma(L+2)\,_{2}\tilde{F}_{3}\left(\frac{L}{2}+1,\frac{L}{2}+\frac{3}{2};L+\frac{5}{2},\frac{b-a}{2}+\frac{L}{2}+\frac{3}{2},\frac{a-b}{2}+\frac{L}{2}+\frac{3}{2};\frac{z}{4}\right)\right)\end{array}\quad.\label{eq:final0_ab-1}\end{equation}
In this form, we must have $a=1,\,3,\,5,\,7\cdots$ and $ $$b=1,\,2,\,3,\cdots,\, a$
. 

One could continue on in this fashion for $p=2$ and beyond, but we will have a product of triplets for  $p=2$ since

\begin{equation}
x^{2}P_{l'}(x)P_{l}(x)=\left(\frac{l'\left(2\left(l'\right)^{2}+l'-3\right)P_{l'-2}(x)+\left(2\left(l'\right)^{3}+5\left(l'\right)^{2}+l'-2\right)P_{l'+2}(x)+\left(4\left(l'\right)^{3}+6\left(l'\right)^{2}-1\right)P_{l'}(x)}{\left(2l'+1\right)\left(4\left(l'\right)^{2}+4l'-3\right)}\right)P_{l}(x)\quad\label{eq:x^2p_lpP_l-1}\end{equation}
and the number of terms increases with $p$. The process of finding the $\,_{2}F_{3}$ products of
sums would, thus, be straight-forward though increasingly complicated as
\emph{p} increases.

\section{Discussion and Concluding Remarks}

We have found that the angular integration of transition amplitudes
for radiative attachment, arising from the Strong Field Approximation
(SFA), provides a means to express certain $\,_{3}F_{4}$ hypergeometric
functions as infinite sums over pair products of $\,_{2}F_{3}$ functions.
The SFA transition amplitudes include products of Generalized Bessel
functions that comprise a product of conventional Bessel
functions, each of which can be expanded in a series of Spherical Harmonics
times $\,_{2}F_{3}$ functions, whose angular integral reduces these
to a single infinite series. On the other hand, one may express the
product of these two conventional Bessel functions as a (different)
$\,_{2}F_{3}$ function, whose angular integral is a $\,_{3}F_{4}$
function. Equating the results of these two methods gives the the
desired relationship.

We have also stepped somewhat away from the physical application of
this relationship by including a multiplicative factor of $\cos^{2}\theta$
to generate a second sort of $\,_{3}F_{4}$ Hypergeometric function
expressed as an infinite sum over pair products of $\,_{2}F_{3}$
functions. If one were to entirely divorce this procedure from transition
amplitudes and replace the Bessel functions with other conventional
functions that can each be expanded in a series of Spherical Harmonics,
and whose product is some $\,_{p}F_{q}$ function, one might be able
to generate additional relations of this general type.

\end{document}